\documentclass[12pt]{amsart}
\usepackage{amsfonts,
amssymb,amsmath,epsfig,graphics,
graphicx,color}
\usepackage{a4wide}
\usepackage[english]{babel}
\newcommand{\halmos}{\rule{1ex}{1.4ex}}
\makeatletter \@addtoreset{equation}{section} \makeatother

\newtheorem{ittheorem}{Theorem}
\newtheorem{itlemma}{Lemma}
\newtheorem{itproposition}{Proposition}
\newtheorem{itcorollary}{Corollary}
\newtheorem{itdefinition}{Definition}
\newtheorem{itremark}{Remark}
\newtheorem{itexamples}{Examples}

\newenvironment{theorem}{\addtocounter{equation}{1}
\begin{ittheorem}}{\end{ittheorem}}
\newenvironment{lemma}{\addtocounter{equation}{1}
\begin{itlemma}}{\end{itlemma}}
\newenvironment{proposition}{\addtocounter{equation}{1}
\begin{itproposition}}{\end{itproposition}}
\newenvironment{corollary}{\addtocounter{equation}{1}
\begin{itcorollary}}{\end{itcorollary}}
\newenvironment{definition}{\addtocounter{equation}{1}
\begin{itdefinition}}{\end{itdefinition}}
\newenvironment{remark}{\addtocounter{equation}{1}
\begin{itremark}}{\end{itremark}}

\newenvironment{examples}{\addtocounter{equation}{1}
\begin{itexamples}}{\end{itexamples}}

\newenvironment{proofs}{\noindent {\em Proof}.\,\,\,}
{\hspace*{\fill}$\halmos$\medskip}

\newenvironment{proofsm}{\noindent}
{\hspace*{\fill}$\halmos$\medskip}

\newcommand{\beq}{\begin{eqnarray}}
\newcommand{\eeq}{\end{eqnarray}}

\newcommand{\beqq}{\begin{eqnarray*}}
\newcommand{\eeqq}{\end{eqnarray*}}

\newcommand{\be}{\begin{equation}}
\newcommand{\ee}{\end{equation}}

\newcommand{\bl}{\begin{lemma}}
\newcommand{\el}{\end{lemma}}

\newcommand{\br}{\begin{remark}}
\newcommand{\er}{\end{remark}}

\newcommand{\bex}{\begin{examples}}
\newcommand{\eex}{\end{examples}}

\newcommand{\bt}{\begin{theorem}}
\newcommand{\et}{\end{theorem}}

\newcommand{\bd}{\begin{definition}}
\newcommand{\ed}{\end{definition}}

\newcommand{\bp}{\begin{proposition}}
\newcommand{\ep}{\end{proposition}}

\newcommand{\bc}{\begin{corollary}}
\newcommand{\ec}{\end{corollary}}

\newcommand{\bpr}{\begin{proofs}}
\newcommand{\epr}{\end{proofs}}

\newcommand{\bprm}{\begin{proofsm}}
\newcommand{\eprm}{\end{proofsm}}

\newcommand{\bi}{\begin{itemize}}
\newcommand{\ei}{\end{itemize}}

\newcommand{\ben}{\begin{enumerate}}
\newcommand{\een}{\end{enumerate}}

\newcommand{\Z}{\mathbb Z}
\newcommand{\R}{\mathbb R}
\newcommand{\N}{\mathbb N}

\newcommand{\T}{\mathbb T}

\newcommand{\cL}{\ensuremath{\mathcal{L}}}

\newcommand{\Ga}{\ensuremath{\Gamma}}

\newcommand{\Om}{\ensuremath{\Omega}}
\newcommand{\calL}{\ensuremath{\mathcal{L}}}

\begin{document}
\title[Couplings and attractiveness II]{Couplings and Attractiveness for 
General Exclusion Processes}
\author{Thierry Gobron}
\address{CNRS UMR 8524, Université de Lille - Laboratoire Paul Painlevé, F-59000 Lille, France}
\author{Ellen Saada}
\address{CNRS, UMR 8145, Laboratoire MAP5\\
Universit\'e Paris Cit\'e, campus Saint-Germain-des-Pr\'es; 
45 rue des Saints-P\`eres
\\75270 Paris cedex 06, France}
\thanks{}
\date{\today}
\keywords{ Particle systems, Attractiveness, Couplings, Discrepancies, 
Invariant measures, Exclusion processes with speed change. }
\subjclass[2000]{Primary 60K35; Secondary 82C22.}

\begin{abstract}
Attractiveness is a fundamental tool to study interacting 
particle systems and the basic coupling 
construction is a usual route to prove this property, as for instance in the 
simple exclusion process. 
 We consider here general exclusion processes where jump rates
from an occupied site to an empty one depend not only on the location of the jump
but also possibly on the whole configuration. These processes include in particular
exclusion processes with speed change introduced by F. Spitzer in 
\cite{spitzer}. For such processes we derive 
necessary and sufficient conditions 
for attractiveness, through the construction 
of a coupled process under which discrepancies do not increase. 
We emphasize the fact that basic coupling is never attractive 
for this class of processes, except in the case  of simple exclusion, 
and that the coupled processes presented here necessarily differ from it. 
We study various examples, for which we determine 
the set of extremal translation invariant
and invariant probability measures. 
\end{abstract}

\maketitle
\setcounter{equation}{0}

\medskip
\hfill{\sl Dedicated to Errico Presutti}
\medskip

\section {Introduction}\label{sec:intro}

Exclusion processes are among the most studied interacting particle systems:
despite their very simple form, these Markov processes exhibit characteristic 
features that make them ideal toy models for many physical or biological phenomena. 

In an exclusion process, particles evolve  on a countable set of sites $S$,  
e.g. $\Z^d$, on which multiple occupancy is forbidden.
 This exclusion rule is encoded in the structure of the
 state space  which is thus  defined as 
 $\Omega=\{0,1\}^S$. For a configuration $\eta\in \Omega$
and for $x\in S$,  $\eta(x)$ is the occupation number at site $x$, that is 
$\eta(x)=1$  whenever a particle is present on site $x$,
while $\eta(x)=0$  when site $x$ is empty.
 Particles jump from one site to another, empty, site 
according to a  probability transition $p(.,.)$ on $S$
(for $S=\Z^d$,  we  consider only translation invariant cases). 

The most widely studied exclusion model is the simple exclusion process (SEP), 
in which particles have all the same speed one, that is the transition rate for 
a particle in a configuration $\eta$ to jump from  its position at site $x$ to 
 an empty site $y$  does not depend on the location of other particles 
 and thus  simply reads $\eta(x)(1-\eta(y))p(x,y)$.
 Endowing  $\Omega$ with the coordinatewise (partial) order, that is,
 for $\eta,\xi\in \Omega$, 
\be\label{eq:order} 
 \eta\leq\xi \Leftrightarrow \forall x\in S,\,\eta(x)\leq\xi(x)
\ee
 we can define  \textit{attractiveness} as the property that this partial order
is maintained through (coupled) evolution whenever it holds at initial time.
 Attractiveness is  a fundamental property of SEP and a  key tool 
 to determine the set $(\mathcal I\cap\mathcal S)_e$ of extremal 
 translation invariant 
 and invariant probability measures for the dynamics
  (see e.g. chapter VIII of \cite{L}). This set consists in 
  a one parameter family $\{\nu_\rho, \rho\in[0,1]\}$
 of Bernoulli product measures, where $\rho$ represents 
 the average particles'density per site. 
 It is also crucial  in establishing hydrodynamics 
 for asymmetric transition probability $p(.,.)$, see e.g. \cite{reza, KL}). 
 In such a problem, attractiveness is  embodied through the 
 ``basic coupling'' construction 
 of two copies  $(\eta_t)_{t\geq 0}$  and $(\xi_t)_{t\geq 0}$
 of simple exclusion processes, under which particles move together 
 as much as possible. 
 In other words, if at some time $s$ particles of both copies 
 attempt to jump, they will try to go from the same departure site $x$ 
 to the same arrival site $y$ according to 
  $p(x,y)$, as long as those jumps are permitted (that is if 
  $\eta_s(x)=\xi_s(x)=1$ and $\eta_s(y)=\xi_s(y)=0$), otherwise only 
  the possible jump will take place. Thanks to basic coupling, it is possible 
  to control the evolution of \textit{discrepancies} between $(\eta_t)_{t\geq 0}$ 
 and $(\xi_t)_{t\geq 0}$, that is, the sites on which the configurations differ.
 Combined with some irreducibility property for
 the probability transition $p(.,.)$, this control is the essential step 
 to derive $(\mathcal I\cap\mathcal S)_e$ (see \cite{L1,L}).
 
 However, ever since the seminal paper \cite{spitzer} by Frank Spitzer in which 
simple exclusion process was first defined, other exclusion processes
have been considered, named exclusion processes with speed change,
 in which jump rates may depend on the configuration
around the particle departure site.
Though such a dependence can be treated within a basic coupling construction 
for (non conservative) spin flip models, it appeared to be not so simple
for conservative ones.
 In order to determine the set $(\mathcal I\cap\mathcal S)_e$ for such models, 
 more involved attractiveness 
 conditions and related coupling constructions were to be found. 
 Sufficient conditions for attractiveness have been obtained by Tom
Liggett in his Saint-Flour lecture notes  \cite{L2} for the models introduced in
\cite{spitzer},
as well as a related coupling leading to $(\mathcal I\cap\mathcal S)_e$ 
whenever these conditions are fulfilled.

 Totally asymmetric versions of exclusion processes  with speed change
  are also natural models of  traffic (see e.g. \cite{GG}). 
Recently, there has been a renewed interest 
in exclusion processes,  in particular those related to integrable models, 
such as the facilitated  exclusion processes
 (see e.g. \cite{BBCS, BESS, AGLS}), or the $q$-Hahn exclusion process 
 (see \cite{BC}).  These models, whether or not attractive, have been analyzed
through other existing techniques such as duality, or through an
 ad-hoc correspondence with (generalized) zero-range processes.
 
  In this work, we consider a general exclusion process on 
   $\Z^d$  and state necessary and sufficient conditions under which 
  attractiveness holds. Here jump rates depend not only on the position 
  and occupation numbers of the sites at which a jump occurs, but also 
  possibly on the whole configuration, so that the basic coupling 
  construction does not hold beyond SEP. 
  We proceed in the spirit of our previous papers 
  on particle systems of misanthrope type \cite{gs, fgs}, in which
  the richer structure of the local state space already imposes 
  non trivial attractiveness conditions even when rates depend 
  on the configuration only through the sites at which a jump occurs. 
  In the course of the present construction, we have to distinguish 
  between the two strongly related notions of monotonicity and attractiveness. 
  Loosely speaking, any two initially ordered configurations which evolve 
  under a monotone process will remain ordered at all time. This property 
  can be stated in a weak sense, or equivalently through the construction 
  of an increasing coupling which preserves the ordering of its marginals. 
  In an attractive process, any pair of configurations behaves in such a way 
  that their differences disappear as much as possible so that they eventually 
  order with probability one. Both notions 
 coincide on classical examples such as the simple exclusion process. We show that
  it is also the case in the  present wider context 
  in the sense that necessary and sufficient conditions for monotonicity imply 
 attractiveness, but additional work has to be done. 
Similarly to simple exclusion,  this property, 
when associated to some irreducibility of 
the coupled process and to an additional assumption on the dynamics, 
eventually leads to a full characterization of extremal translation invariant,
 invariant probability measures of generalized exclusion processes.

 The paper is organized as follows. In Section \ref{sec:model-results} 
we define the generalized exclusion model, and state our main results:
necessary and sufficient conditions for monotonicity (Theorem \ref{thm:1}), 
the existence, for a monotone process, of an increasing coupling 
under which discrepancies do not increase (Theorem \ref{thm:2}),
 and  determination of the set $(\mathcal I\cap\mathcal S)_e$ 
(Theorem \ref{prop:inv-meas}). 
In Section \ref{sec:couplings}, we prove Theorem \ref{thm:1} and give 
in a series of propositions the
construction of the successive generators leading to Theorem \ref{thm:2}
 and Theorem \ref{prop:inv-meas}). 
These propositions as well as Theorem \ref{thm:2}
are proved in Section \ref{sec:proofs}.
In Section \ref{sec:Applications},  we illustrate our results with examples, 
showing first that our construction reduces to 
basic coupling in the case of simple exclusion and only there.
We then consider exclusion processes with speed change, extending the results of  
\cite{spitzer, L2}. 
Finally, we turn to traffic models, considering first a generalization of
 the totally asymmetric 2-step exclusion process studied in \cite{guiol}, 
 and a symmetrized version of the totally asymmetric traffic model from \cite{GG}. 
In all cases, we compute explicitely the  attractive coupling rates  
and give the set of invariant measures $(\mathcal I\cap\mathcal S)_e$.

\section{Model and Main Results}\label{sec:model-results}
\setcounter{equation}{0}
In this section,  we  define the class of exclusion models we consider and 
state our two main results:
Theorem \ref{thm:1} gives necessary and sufficient conditions 
for monotonicity, and Theorem \ref{thm:2} links monotonicity and 
attractiveness for this model, through a coupling construction.

We first introduce a general exclusion process $(\eta_t)_{t\ge 0}$ on $S=\Z^d$, 
together with some notation and general properties. 
Let $\Omega=\{0,1\}^S$ be
its state space and $\cL$ its formal generator, acting on any cylinder 
function $f$ and for any configuration $\eta\in\Omega$,
\begin{eqnarray}
\label{generator}
\cL f(\eta)&=&\sum_{x,y \in S}\eta(x)(1-\eta(y))\Ga_\eta(x,y)\bigl[f(\eta^{x,y}) - f(\eta)\bigr]
\end{eqnarray}
where for any $(x,y)\in S^2$, $\eta^{x,y}$ is a the configuration 
obtained from $\eta$ by exchanging the occupation numbers in configuration
 $\eta$ at  sites $x$ and $y$
\be
\eta^{x,y}(z) =
\begin{cases}
\eta(y) &\hbox{ if }  z=x\cr
\eta(x) &\hbox{ if }  z=y\cr
\eta(z) &\hbox{ otherwise}\cr
\end{cases}
\ee
 The process is thus  conservative, and the quantity $\eta(x)+\eta(y) $ 
is conserved in a jump from site $x$ to site $y$. We denote by $(T(t),t\geq 0)$
the semi-group of this process. 
 \begin{remark}\label{rk:SEP}
When the jump rates $\Ga_\eta(x,y)$ are independent of the configuration
 $\eta$, and reduce to a probability transition $(p(x,y),x,y\in S)$ on $S$, 
\be\label{rate:SEP}
\Ga_\eta(x,y)=p(x,y)
\ee
one recovers the simple exclusion process.
\end{remark} 
We assume the following conditions on the jump rates, so that
 \eqref{generator} is the infinitesimal generator 
of a well defined Markov process (see \cite[Chapter I]{L}):
\begin{equation}\label{existence}
\sup_{v\in S}\sum_{u\in S}\sup_{\eta\in\Omega}\Ga_\eta(u,v)<+\infty 
\quad\hbox{ and }\quad
\sup_{u\in S}\sum_{v\in S}\sup_{\eta\in\Omega}\Ga_\eta(u,v)<+\infty
\end{equation}  
 Of course, these generic conditions can be alleviated, 
 depending on the example at hand.

Let us recall the \textit{monotonicity property} for particle systems, quoting 
\cite[Chapter II]{L}. We denote by $\mathcal M$ the set of all bounded,
non-decreasing, continuous 
functions $f$ on $\Omega$. The partial order \eqref{eq:order} induces 
a stochastic order on the set $\mathcal P$ of probability measures on $\Omega$
endowed with the weak topology:
\be\label{eq:order-on-P}
\forall \nu,\nu'\in{\mathcal P},\,
\nu\leq\nu'\Leftrightarrow\bigl(\forall f\in{\mathcal M},\nu(f)\leq\nu'(f)\bigr)
\ee
\bt\label{th:lig-attra} \cite[Chapter II, Theorem 2.2]{L}
For the particle system $(\eta_t)_{t\geq 0}$ the following 
two statements are equivalent.

(a) $f\in\mathcal M$ implies $T(t)f\in\mathcal M$ for all $t\geq 0$. 

(b)  For $\nu,\nu'\in{\mathcal P}$, $\nu\leq\nu'$ implies $\nu T(t)\leq\nu'T(t)$
for all $t\geq 0$.
\et 
\bd\label{def:monotonicity} \cite[Chapter II, Definition 2.3]{L}
The particle system  $(\eta_t)_{t\ge 0}$ is {\rm monotone}
 if the equivalent statements of Theorem \ref{th:lig-attra} are satisfied.
\ed
 Our first 
 main result is the following set of necessary and sufficient 
conditions for monotonicity.
\bt\label{thm:1}
The exclusion process defined by 
\eqref{generator} is {\rm monotone} if and only if 
for any couple of configurations 
$(\xi,\zeta)\in\Omega^2$ such that $\xi\le\zeta$,
the following inequalities hold:\\
For all $y\in S$ such that $\zeta(y)=0$,
       \begin{equation}
       \label{eq:1}
 \sum_{x\in S} \xi(x) \bigl[\Ga_\xi(x,y) - \Ga_\zeta(x,y)\bigr]^+ 
\le
       \sum_{x\in S} \zeta(x) (1 - \xi(x) ) \Ga_\zeta(x,y)  
       \end{equation}
       For all $x\in S$ such that $\xi(x)=1$,
\begin{equation}
       \label{eq:2}
\sum_{y\in S} (1 - \zeta(y) )\bigl [\Ga_\zeta(x,y) - \Ga_\xi(x,y)\bigr]^+  
\le
       \sum_{y\in S} \zeta(y) (1 - \xi(y) ) \Ga_\xi(x,y)       
     \end{equation}
\et
In  Section \ref{sec:couplings}, we prove that these conditions 
are necessary and rely on them to build in Propositions \ref{propcoupling0}
and \ref{propincreasing} 
a coupling between two copies of the process. 
A coupling is called \textit{increasing} if it preserves 
the stochastic order  between marginal  configurations. 
In Section \ref{sec:proofs},  we achieve the proof of Proposition
\ref{propincreasing}, that is,  this coupling is proven to be \textit{increasing}
under the hypothesis that inequalities \eqref{eq:1}--\eqref{eq:2} hold, 
showing in turn that these conditions are also sufficient.

Beyond monotonicity, a coupling construction turns out 
to be essential to characterize the set 
$\left(\mathcal {I}\cap\mathcal{S}\right)_e$ of 
extremal invariant and translation invariant probability measures
of $(\eta_t)_{t\ge 0}$.  In our setting, the marginals of the coupled process 
built in Propositions \ref{propcoupling0}
and \ref{propincreasing}  are not necessarily ordered, and
the evolution of the discrepancies  between them is the main object to control:
\bd\label{def-disc}
In a  coupled process  $(\xi_t,\zeta_t)_{t\ge 0}$, there is a {\rm discrepancy} 
 at site $z\in S$ at time $t$ if $\xi_t(z)\ne\zeta_t(z)$. 
 \ed
A process is  deemed
{\sl attractive} if there exists a coupling of 
two copies of the process such that the discrepancies 
between the marginals do not increase in time 
and eventually disappear, so that the marginals will eventually become ordered 
 with probability one. 

 For instance, as recalled in 
the introduction, 
basic coupling is attractive for the simple exclusion process (SEP)
 and in any coupled transition the number of discrepancies  on the involved sites 
 remains constant whenever the values of the two marginal 
 configurations are  ordered, but decreases otherwise. 

Beyond this case, an increasing coupling does not necessarily 
impose constraints on the coupled evolution of 
unordered pairs of configurations, so that 
the number of discrepancies is not necessarily  non-increasing. 
However here we have the following:
\bt\label{thm:2}
Suppose that the process defined by 
\eqref{generator} is  monotone  on $\Omega=\{0,1\}^S$. 
Then  it is attractive, that is, there exists an increasing  coupled process
on $\Omega\times\Omega$ such that the number of
discrepancies does not increase with time.
\et
The proof of Theorem \ref{thm:2} relies on the explicit construction of 
such an {\sl attractive} coupling,  which refines the previous 
{\sl increasing} one. 
It is described in Proposition \ref{propdecdisc}, 
while proofs of existence and attractiveness are postponed 
to Section \ref{sec:proofs}.

Therefore in our setting, monotonicity and attractiveness coincide, 
so that we will speak only of attractiveness when dealing with examples 
in Section \ref{sec:Applications}. 

To conclude with the characterization of the set 
$\left(\mathcal {I}\cap\mathcal{S}\right)_e$, we need not only that
in the coupling process the number of
discrepancies does not increase with time, but also that this number
decreases. For this, we need to construct again another coupling process.
But it requires an additional assumption of the dynamics.
\bd\label{def:blocking-conf}
 An exclusion process with generator \eqref{generator} 
 has no {\rm blocking configurations} if for  any
 configuration $\xi\in\Omega$, the set of {\rm open edges} 
 $\{(x,y)\in S^2 : \Gamma_\xi(x,y)>0 \}$ is independent on $\xi$. 
The set $S$ is then said {\rm fully connected } if  
for all $(x,y)\in S^2$,$x\not= y$, there exists a finite open path in $S$
between $x$ and $y$, that is a sequence $\{x_0,\cdots,x_n\}$ 
for some $n>0$ such that $(x_{i-1},x_i)$ is open for 
$i\in\{1,\cdots,n\}$ with either  $x_0=x$ and $x_n=y$, 
or $x_0=y$ and $x_n=x$.
\ed
In Subsection \ref{subsec:inv-meas}, we will explain how,
whenever the dynamics has
no blocking configurations and $S$ is fully connected, 
it is possible to construct a coupling such that any pair of discrepancies of opposite 
sign have a positive probability to disappear in finite time.
When the jump rates are translation invariant, this reduces the 
derivation of the set $\left(\mathcal {I}\cap\mathcal{S}\right)_e$ 
essentially to the classical proof, originally applied to the
simple exclusion process (going back to \cite{L1}),
which leads to the following theorem. 
 \bt\label{prop:inv-meas}
Let $(\eta_t)_{t\ge 0}$ be an exclusion process 
with generator \eqref{generator} and translation invariant jump rates,
such that there are no blocking configurations
and $S$ is fully connected in the sense of Definition \ref{def:blocking-conf}.
If $(\eta_t)_{t\ge 0}$
is attractive then 

 1) The set of translation invariant,  extremal invariant measures  
 $(\mathcal I\cap\mathcal S)_e$ is a one parameter family 
 $ \{\mu_\rho,\,\rho\in{\mathcal R}\}$, 
where $\mathcal R$ is a closed subset of $[0,1]$ containing 
$\{0,1\}$, and for every $\rho \in\mathcal R$,  $\mu_\rho$ is 
a translation invariant probability 
measure on $\Omega$ with $\mu_\rho[\eta(0)]=\rho$; furthermore, the measures 
$\mu_\rho$ are stochastically ordered, that is,
$\mu_\rho\leq\mu_{\rho'}$ if $\rho\leq\rho'$;

 2)  if $(\eta_t)_{t\ge 0}$ possesses a one parameter family
$\{\mu_\rho\}_\rho$ of product invariant and translation 
invariant probability measures, we have 
$(\mathcal I\cap\mathcal S)_e=\{\mu_\rho\}_\rho$.
 \et

Our  results can be extended in various ways, to more general conservative models, 
as well as to some mixed non conservative models with both exchanges and configuration 
independent birth-death events, but this is beyond the scope of the present paper. 


\section{Proofs of main theorems and coupling constructions.}\label{sec:couplings}
 This section is devoted to the construction of the coupling necessary to 
the proof of Theorem \ref{thm:1}
(in Subsection \ref{subsec:proof-thm1}),
in three steps. 
We first prove that inequalities \eqref{eq:1}--\eqref{eq:2}
are necessary  conditions. In order to prove that these conditions 
 are also sufficient, we introduce  in Proposition \ref{propcoupling0}
 the general form $\bar\cL$  of a Markovian coupling generator 
 associated to $\cL$, depending on a set of coupled transition rates 
 $G_{\xi,\zeta}(.)$. Those rates are defined in Proposition \ref{propincreasing} 
and we prove in turn that with such a choice, and whenever inequalities 
\eqref{eq:1}--\eqref{eq:2} 
are fulfilled, the generator $\bar\cL$ defines an increasing coupling.
We continue this section (in Subsection \ref{subsec:proof-thm2})
with the proof of Theorem \ref{thm:2}, introducing in Proposition 
\ref{propdecdisc} the generator ${\overline \calL}^D$
 of an attractive coupling. Finally we explain in Subsection \ref{subsec:inv-meas}
how to prove Theorem \ref{prop:inv-meas} by refining the construction of an attractive
 coupling (in Proposition \ref{decrdiscr}). 
 Proofs of the above Propositions
 are given in Section \ref{sec:proofs}. 

\subsection{Proof of Theorem \ref{thm:1}}\label{subsec:proof-thm1}
 Inequalities \eqref{eq:1}--\eqref{eq:2} are particular instances 
 (and in turn the worst cases) 
of a larger set of inequalities  (first derived by A.W. Massey \cite{M}) 
that the coefficients of the infinitesimal generator of a monotone Markov 
process need to fulfil. 
We sketch their derivation hereafter  and we refer to \cite{M} 
for a thorough derivation (see also \cite{gs} for details). 
The idea is to derive sensible necessary conditions on the jump
rates for a Markov process to be monotone, using the fact that 
the characteristic function of any increasing (or decreasing) cylinder  
set $V\subset \Omega$, is a monotone cylinder function 
on $\Omega$. Let  $(\xi_t)_{t\geq 0}$ and $(\zeta_t)_{t\geq 0}$ 
two instances of a monotone process 
with initial conditions $\xi_0$ and $\zeta_0$ such that $\xi_0 \le \zeta_0$, 
then ${\bf 1}_V(\xi_t) \le {\bf 1}_V(\zeta_t)$  (and reverse inequality 
for a decreasing set). In addition if initial conditions are chosen so that 
 $\xi_0 \not\in V$ and $ \zeta_0 \not\in V$, the same inequality 
 holds for the ratios
\[
\frac{1}{t} \left( {\bf 1}_V(\xi_t)- {\bf 1}_V(\xi_0)\right)
 \le\frac{1}{t} \left( {\bf 1}_V(\zeta_t)- {\bf 1}_V(\zeta_0)\right)
\]
for all $t>0$. Taking properly the limit $t\to 0$ gives then 
inequalities involving the rates 
of the Markov generator, hereafter named  
``Massey conditions'' and stated below in our case:

{\sl If the particle system defined in \eqref{generator} is monotone, then
for all configurations $(\xi,\zeta)\in\Om\times\Om$ such
that $\xi\le\zeta$,
\ben
\item For all increasing cylinder sets $V\subset \Om$ such that
$\zeta\notin V$,
\be\label{massey-incr}
\sum_{x,y} \xi(x)(1-\xi(y)) \Ga_\xi(x,y)
1_V(\xi^{x,y})  
\leq
\sum_{x,y} \zeta(x)(1-\zeta(y)) \Ga_\zeta(x,y)
1_V(\zeta^{x,y})
\ee
\item For all decreasing cylinder sets $V\subset \Om$ such that
$\xi\notin V$,
\be\label{massey-decr}
\sum_{x,y}   \zeta(x)(1-\zeta(y)) \Ga_\zeta (x,y) 1_V(\zeta^{x,y})
\leq
\sum_{x,y} \xi(x)(1-\xi(y))  \Ga_\xi (x,y) 1_V(\xi^{x,y})
\ee
\een }
\mbox{}\\ 
\bpr [Theorem \ref{thm:1}, Necessary conditions].
\mbox{}\\ 
Equations \eqref{eq:1} follow from \eqref{massey-incr} by taking 
a particular sequence of cylinder
increasing sets and passing to the limit. Equations \eqref{eq:2} 
follow in the same way from \eqref{massey-decr}.
Let $\xi$, $\zeta$ be two configurations such that $\xi\le \zeta$ 
and take $y$ such that $\zeta(y)=0$.
For $n>0$, we construct a configuration $\eta_n$ as follows.
\begin{eqnarray}
\label{def:etan}
\eta_n(x)=
  \begin{cases}
    1 & \text{ if } x=y, \\
    1 & \text{ if }  \|x-y\|\le n\text{ , } \xi(x)=1 
    \text{ and } \Ga_\xi(x,y)<\Ga_\zeta(x,y), \\
    0 & \text{otherwise}.
  \end{cases}
\end{eqnarray}
We define the increasing cylinder set $V_n=\{\rho\in\Om, \rho\ge\eta_n\}$. 
 Since $\zeta(y)=0$, configuration $\zeta$ (and hence $\xi$) does not belong 
 to $V_n$.  Equation \eqref{massey-incr} applied to $V_n$ now selects single 
 jumps which allow to enter $V_n$, hence moving a particle from any site $x$ 
 with $\eta_n(x)=0$ to site $y$. We thus get:
\begin{eqnarray}
\label{incr2}
\sum_{x\in S}
\xi(x) (1-\eta_n(x)) \Ga_\xi(x,y)
\leq
\sum_{x\in S}
\zeta(x) (1-\eta_n(x)) \Ga_\zeta(x,y)
\end{eqnarray}
 Note that by conditions \eqref{existence},  both sums are finite.
For all $x\ne y$,  we have
\begin{eqnarray*}
\zeta(x) (1-\eta_n(x))&=& \zeta(x) (1-\eta_n(x)) (1-\xi(x)) 
+  \zeta(x) (1-\eta_n(x)) \xi(x) \nonumber\\
&=& \zeta(x) (1-\xi(x)) + \xi(x) (1-\eta_n(x)) 
\end{eqnarray*}
where the second line comes from the fact that $\eta_n(x)\le \xi(x) \le\zeta(x) $.
Inserting this expression in the right hand side of \eqref{incr2}, we get
\begin{eqnarray*}
\sum_{x\in S}
 \xi(x) (1-\eta_n(x))  \bigl(\Ga_\xi(x,y)-\Ga_\zeta(x,y)\bigr)
\leq
\sum_{x\in S}
\zeta(x) (1-\xi(x)) \Ga_\zeta(x,y)
\end{eqnarray*}
which gives, using definition \eqref{def:etan} of $\eta_n$
\beq
\label{eq:10}
\sum_{x\in S : |x-y|\le n}
\xi(x) \bigl[\Ga_\xi(x,y)-\Ga_\zeta(x,y)\bigr]^+ +
\sum_{x\in S : |x-y|> n}
\xi(x)  \bigl(\Ga_\xi(x,y)-\Ga_\zeta(x,y)\bigr)\nonumber\\
\qquad\qquad\leq
\sum_{x\in S}
\zeta(x) (1-\xi(x)) \Ga_\zeta(x,y)
\eeq
 Conditions \eqref{existence} now imply that the second term in the left hand side 
 of \eqref{eq:10} goes to zero as $n\to\infty$.
 Taking the limit $n\to\infty$ in \eqref{eq:10} thus gives 
\[
\sum_{x\in S}
\xi(x) \bigl[\Ga_\xi(x,y)-\Ga_\zeta(x,y)\bigr]^+ \leq
\sum_{x\in S}
\zeta(x) (1-\xi(x)) \Ga_\zeta(x,y)
\]
which is Equation \eqref{eq:1}.

 Equation \eqref{eq:2} can be derived in a similar way from \eqref{massey-decr}. 
Let again $\xi$, $\zeta$ be two configurations such that $\xi\le \zeta$ 
and take now $x\in S$ such that $\xi(x)=1$.
Let $n>0 $ and consider the configuration
$\eta_n$ such that:
\begin{eqnarray}
\label{def:etanbis}
\eta_n(x)=
  \begin{cases}
    0 & \text{ if } y=x\, \\
        0 & \text{ if } |x-y|\le n\text{ , }\zeta(y)=0 
        \text{ and }\Ga_\xi(x,y)>\Ga_\zeta(x,y), \\
    1 & \text{otherwise}.
  \end{cases}
\end{eqnarray}
We construct the decreasing cylinder set $V_n=\{\rho\in\Om, \rho\le\eta_n\}$. 
Since $\xi(x)=1$, the configuration $\xi$ (and thus $\zeta$) does not 
belong to $V_n$.
Equation \eqref{massey-decr} now selects single jumps which allow 
to enter the decreasing set, thus removing a particle at $x$ and moving it
 to any possible site $y$ where $\eta_n(y)=1$. We thus get
\be\label{decr2}
\sum_{y\in S} \eta_n(y) (1-\zeta(y)) \Ga_\zeta (x,y)
\leq
\sum_{y\in S} \eta_n(y) (1-\xi(y))  \Ga_\xi (x,y)
\ee
For all $y\ne x$,  we now have
\begin{eqnarray*}
\eta_n(y) (1-\xi(y) )&=&\eta_n (y) (1-\xi(y) )(1-\zeta(y))
 + \eta_n (y)(1-\xi(y)) \zeta (y)\\
&=& \eta_n (y)(1-\zeta(y)) + \zeta (y)(1-\xi(y))\nonumber
\end{eqnarray*}
where we have used that $ \eta_n (y) \ge  \zeta (y) \ge \xi(y) $.
Inserting this expression in the right hand side of \eqref{decr2} gives
\[
\sum_{y\in S} \eta_n(y) (1-\zeta(y))\bigl( \Ga_\zeta (x,y))-\Ga_\xi(x,y)\bigr)
\leq
\sum_{y\in S} \zeta(y)  (1-\xi(y))  \Ga_\xi (x,y)
\]
Using the definition  \eqref{def:etanbis}  of $\eta_n$, we get
\beq
\label{eq:11}
\sum_{y\in S: |y-x|\le n} (1-\zeta(y)) \bigl[\Ga_\zeta (x,y) - \Ga_\xi (x,y)\bigr]^+
+
\sum_{y\in S: |y-x|> n}(1-\zeta(y)) \bigl(\Ga_\zeta (x,y) 
- \Ga_\xi (x,y)\bigr)\nonumber\\
\qquad\qquad\leq
\sum_{y\in S}\zeta(y) (1-\xi(y) ) \Ga_\xi (x,y)
\eeq
In the limit $n\to\infty$, the second term in the left hand side of \eqref{eq:11}
goes to zero and one gets
\beq
\nonumber
\sum_{y\in S} (1-\zeta(y)) \bigl[\Ga_\zeta (x,y) - \Ga_\xi (x,y)\bigr]^+
\leq
\sum_{y\in S}\zeta(y) (1-\xi(y) ) \Ga_\xi (x,y)
\eeq
which is Equation \eqref{eq:2}.
 \epr
 
Equations \eqref{eq:1}--\eqref{eq:2} can be interpreted in the following way. 
First,  by conditions \eqref{existence}, the sums 
appearing in \eqref{eq:1}--\eqref{eq:2} are always finite. The right hand side 
of \eqref{eq:1} measures the excess rate at which an empty site $y$ is filled 
in the smaller configuration $\xi$, so that  coupling jumps in both 
configurations from the same initial sites $x$ to $y$ will be clearly not 
sufficient to preserve partial order if this sum is different from zero. 
Equation \eqref{eq:1} suggests that partial order could be preserved by 
coupling such ``excess rate" jumps with jumps involved in the left hand side, 
that is jumps to $y$ from sites occupied in configuration $\zeta$, but empty 
in $\xi$. Equation \eqref{eq:1} just states that such rates are sufficient to do so. 

Equation \eqref{eq:2}  can be interpreted in a similar way: Now the right hand 
side measures the excess rate at which a filled site $x$ is depleted in the 
larger configuration $\zeta$, so that again partial order could not be preserved 
by coupling jumps in both configurations from site $x$ to the same site $y$ 
whenever this sum differs from zero. Again equation \eqref{eq:2} suggests that 
partial order could be preserved by coupling this second set of ``excess 
rate jumps'' with jumps in the smaller configuration $\xi$ from the same site 
$x$ to any site $y$, empty in configuration $\xi$ but already filled in $\zeta$. 
Again equation \eqref{eq:2} states that the jump rates are just sufficient to do so.

  We now use these ideas to construct a  coupling process 
  then prove that it is increasing, that is, 
  we proceed with the second and third steps of the proof of  Theorem \ref{thm:1}. 
  
\bpr [Theorem \ref{thm:1}, Coupling Process].
\mbox{}\\ 
We define the general form an increasing coupling process should take. 
\bp\label{propcoupling0}
The operator  $\bar\cL$ defined, for any cylinder function $f$ on $\Om\times\Om$
 and any pair of configurations  $(\xi,\zeta)\in\Om\times\Om$, by 
\beq\label{coupling}
\overline{\calL}f(\xi,\zeta)&=&
\sum_{x_1,y_1 \in S}  \xi(x_1)  (1-\xi(y_1)) \Ga_\xi(x_1,y_1)
\bigl(f(\xi^{x_1,y_1},\zeta) - f(\xi,\zeta)\bigr)\nonumber\\
&+&\sum_{x_2,y_2 \in S}  \zeta(x_2)  (1-\zeta(y_2)) \Ga_\zeta(x_2,y_2)
\bigl(f(\xi,\zeta^{x_2,y_2}) - f(\xi,\zeta)\bigr)\\
&+&\sum_{x_1,y_1\in S}  \sum_{x_2,y_2\in S} 
 \xi(x_1)  (1-\xi(y_1))  \zeta(x_2)  (1-\zeta(y_2)) 
G_{\xi,\zeta}(x_1,y_1;x_2,y_2)\nonumber\\
&&\qquad\qquad\qquad\times 
\bigl(f(\xi^{x_1,y_1},\zeta^{x_2,y_2}) - f(\xi^{x_1,y_1},\zeta) 
- f(\xi,\zeta^{x_2,y_2}) +f(\xi,\zeta)\bigr)\nonumber
\eeq
is the generator of a Markovian coupling between two copies of the Markov process 
defined by \eqref{generator}, provided that
 for all pairs of configurations $(\xi,\zeta)\in\Omega^2$ the coupling rates  
 $ G_{\xi,\zeta}$ are non-negative and  the following inequalities hold 
\beq\label{ineqcoupling1}
\forall (x_1,y_1)\in S^2,  
\sum_{x_2,y_2\in S}  \zeta(x_2)  (1-\zeta(y_2))\, G _{\xi,\zeta}(x_1,y_1;x_2,y_2)
\le \Ga_\xi(x_1,y_1) \\
\forall (x_2,y_2)\in S^2,  
\sum_{x_1,y_1\in S}  \xi(x_1)  (1-\xi(y_1))\,  G_{\xi,\zeta}(x_1,y_1;x_2,y_2)
\le \Ga_\zeta(x_2,y_2)\label{ineqcoupling2}
\eeq
\ep
Proof of Proposition \ref{propcoupling0} is postponed to Section 
\ref{sec:proofs}. As a shorthand notations for the sums appearing in 
the left hand side of equations \eqref{ineqcoupling1}--\eqref{ineqcoupling2},
 we define for all couples of configurations $(\xi,\zeta)\in \Om\times\Om$ 
 and  all $(x,y)\in S^2$, the quantities 
\beq
\label{coudis3}
\varphi_{\xi,\zeta} (x,y) &:=&
\sum_{x',y'\in S}  \zeta(x')  (1-\zeta(y'))   G_{\xi,\zeta} (x,y;x',y')
\\
\label{coudis2}
\overline\varphi_{\xi,\zeta} (x,y) &:=& 
\sum_{x',y'\in S}  \xi(x')  (1-\xi(y'))   G_{\xi,\zeta} (x',y';x,y)
\eeq 
\epr

\bpr  [Theorem \ref{thm:1}, Increasing Coupling Process].
\mbox{}\\ 
We now give the set of coupling rates $G_{\xi,\zeta} (x,y;x',y')$ 
which defines an increasing coupling.

We first introduce some notations.
Let $\xi$ and $\zeta$ be two configurations in $\Om$. 
For all $x\in S$ such that $\xi(x) =\zeta(x)=1$,
we define the two sets
\beq\label{Yx}
Y_{\xi,\zeta}^x&=&\{y\in S : \xi(y)=0 , \zeta(y)=1, \Gamma_\xi(x,y) >0 \}
\\\label{bYx}
{\overline Y}_{\xi,\zeta}^x&=&\{y\in S : \xi(y)=\zeta(y)=0, 
\Gamma_\zeta(x,y) >\Gamma_\xi(x,y) \}
\eeq
Whenever they are non empty, we define an arbitrary order on these two sets, 
possibly depending on $\xi$, $\zeta$ and $x$, and denote by 
$y_{\xi,\zeta}^{x,k}$ (respectively  $\overline y_{\xi,\zeta}^{x,k}$) 
the $k^\text{\rm th}$ element in $Y_{\xi,\zeta}^x$ (respectively 
$\overline Y_{\xi,\zeta}^x$).

Similarly, for all $y\in S$ such that $\xi(y) =\zeta(y)=0$, we define
\beq\label{Xy}
X_{\xi,\zeta}^y&=&\{x\in S : \xi(x)=\zeta(x)=1, \Gamma_\xi(x,y) 
>\Gamma_\zeta(x,y) \}
\\\label{bXy}
\overline X_{\xi,\zeta}^y&=&\{x\in S : \xi(x)=0 , \zeta(x)=1, 
\Gamma_\zeta(x,y) >0 \}
\eeq
We define an arbitrary order on these two sets as well, 
possibly depending on  $\xi$, $\zeta$ and $y$, and denote
by $x_{\xi,\zeta}^{y,k}$ (respectively  $\overline x_{\xi,\zeta}^{y,k}$) 
the $k^\text{\rm th}$ element in $X_{\xi,\zeta}^y$ (respectively 
$\overline X_{\xi,\zeta}^y$). 

For definiteness, when one of the above sets is finite or empty,  
say $|Y_{\xi,\zeta}^x|=C_Y<\infty$, we may extend the ordered sequence 
of its elements to an infinite one, $(y_{\xi,\zeta}^{x,n})_{n>0}$, by 
setting arbitrarily $ y_{\xi,\zeta}^{x,n}= 0$ for all $n> C_Y$.

For all $x\in S$ such that $\xi(x) =\zeta(x)=1$, we define the two 
series $\left(S_{\xi,\zeta}^{x,n}\right)_{n\ge 0} $
and $\left(\overline T_{\xi,\zeta}^{x,n}\right)_{n\ge 0} $ such that 
$S_{\xi,\zeta}^{x,0} = 0$, $\overline T_{\xi,\zeta}^{x,0} = 0$, and
\beq\label{SYx}
S_{\xi,\zeta}^{x,n}&=&
 \sum_{k=1}^{n\wedge |Y_{\xi,\zeta}^x|} \Gamma_\xi (x,y_{\xi,\zeta}^{x,k}) 
 \qquad \forall n>0
\\
\overline T_{\xi,\zeta}^{x,n} &=&
 \sum_{k=1}^{n\wedge |\overline Y_{\xi,\zeta}^x|} 
 \left[\Gamma_\zeta (x, \overline y_{\xi,\zeta}^{x,k}) 
 - \Gamma_\xi (x, \overline y_{\xi,\zeta}^{x,k}) \right]^+
 \qquad \forall n>0\label{bTYx}
\eeq
Similarly, for all $y\in S$ such that $\xi(y) =\zeta(y)=0$, we define 
the two series $\left(T_{\xi,\zeta}^{y,n}\right)_{n\ge 0} $
and $\left(\overline S_{\xi,\zeta}^{y,n}\right)_{n\ge 0} $ such that 
$T_{\xi,\zeta}^{y,0} = \overline S_{\xi,\zeta}^{y,0} = 0$ and
\beq\label{TXy}
T_{\xi,\zeta}^{y,n} &=&
 \sum_{k=1}^{n\wedge |X_{\xi,\zeta}^y|}
\left[ \Gamma_\xi(x_{\xi,\zeta}^{y,k},y) 
-\Gamma_\zeta(x_{\xi,\zeta}^{y,k},y)\right]^+ \qquad \forall n>0
\\
\overline S_{\xi,\zeta}^{y,n} &=&
 \sum_{k=1}^{n\wedge |\overline X_{\xi,\zeta}^y|}
 \Gamma_\zeta(\overline x_{\xi,\zeta}^{y,k},y) \qquad \forall n>0
\label{bSXy}
\eeq
Note that by definition,   the four series have nonnegative terms 
and are nondecreasing, and by
\eqref{existence}, they are also convergent.

Finally, for any two  convergent 
series $(S_n)_{n\ge 0}$ and   $(T_n)_{n\ge 0}$, 
we define the quantity $H_{n,m}(S,T)$
for all $n >0$ and all $m>0$ as
\be\label{Hmn}
H_{m,n} (S_.,T_.) = S_m \wedge T_n - S_{m-1} \wedge T_n - S_m \wedge T_{n-1}
 + S_{m-1} \wedge T_{n-1}
\ee
Note that $H_{m,n} (S_.,T_.)\ge 0$ whenever $S_.$ and $T_.$ 
are nondecreasing series.

We can now state the following
\bp\label{propincreasing}
Under conditions \eqref{eq:1}--\eqref{eq:2}, 
the generator  given by  \eqref{coupling} with coupling rates  
$ G_{\xi,\zeta}$ below, 
defines an increasing Markovian coupling.
\beq\label{G-rates}
&&G _{\xi,\zeta}(x,y;x',y') =\nonumber\\
&&\qquad\qquad\begin{cases}
 \delta(x,x')\, \delta(y,y')\, \Ga_\xi(x,y)\wedge\Ga_\zeta(x,y)\\
\qquad+  \delta(x,x'){\displaystyle \sum_{m,n>0} 
\delta(y,y_{\xi,\zeta}^{x,m}) \, \delta(y',\overline y_{\xi,\zeta}^{x,n})} \,
 H_{m,n} (S_{\xi,\zeta}^{x,.},\overline T_{\xi,\zeta}^{x,.})\\
\qquad+ \delta(y,y') {\displaystyle\sum_{m,n>0}
 \delta(x,x_{\xi,\zeta}^{y,m}) \delta(x',\overline x_{\xi,\zeta}^{y,n}) }\,
H_{m,n} (T_{\xi,\zeta}^{y,.},\overline S_{\xi,\zeta}^{x,.}) &\text{ if } \xi\le\zeta
\cr
 \delta(x,x') \,\delta(y,y')\, \Ga_\xi(x,y)\wedge\Ga_\zeta(x,y)\\
\qquad+  \delta(x,x'){\displaystyle \sum_{m,n>0} }
  \delta(y,\overline y_{\zeta,\xi}^{x,m}) \, \delta(y', y_{\zeta,\xi}^{x,n}) \, 
H_{m,n} (\overline T_{\zeta,\xi}^{x,.},S_{\zeta,\xi}^{x,.})\\
\qquad+ \delta(y,y') {\displaystyle\sum_{m,n>0}}
 \delta(x,\overline x_{\zeta,\xi}^{y,m})\, \delta(x',x_{\zeta,\xi}^{y,n}) \,
H_{m,n} (\overline S_{\zeta,\xi}^{y,.},T_{\zeta,\xi}^{y,.}) &\text{ if } \xi>\zeta
\cr
0  &\text{ otherwise }
\end{cases}
\eeq
\ep
\br\label{uncoupled}  
With the above choice, jumps are uncoupled
unless $\xi$ and $\zeta$ are ordered. In such a case,
the coupling rate $G _{\xi,\zeta}(x,y;x',y')$ is possibly non zero 
only if the two coupled jumps have either the same
initial point,  the same final point, or both.
\er
\br\label{order}  
The ordering in the four ensembles defined in \eqref{Yx}--\eqref{bXy} 
can be chosen arbitrarily, possibly as a function of the configurations $\xi$ 
and $\zeta$ and on the (initial or final) common jump site. The best choice
may depend on the particular system at hand, and different choices lead 
to different increasing couplings.
Furthermore, all these coupling are extremal in the sense that they cannot 
be written as a convex combination of   
other increasing couplings, while any convex combination of these is 
again an increasing coupling.
\er
\br\label{G=0}
In definition \eqref{G-rates}, the first (resp. second) sum appearing in the 
right hand side in the case $\xi\le\zeta$
is zero except possibly when there is a jump in the first marginal $\xi$ 
from a site $x$ to a site $y\in Y_{\xi,\zeta}^x$  
coupled with a jump  in the second marginal $\zeta$ from the same site $x$ 
to a site $y'\in \overline Y_{\xi,\zeta}^x$
(respectively a jump in the first marginal from a site in $X_{\xi,\zeta}^y$ 
coupled to a jump in the second marginal
from a site in $\overline X_{\xi,\zeta}^y$ to the same site $y$).
Moreover, by the definitions \eqref{Yx}--\eqref{bYx} of  
$Y_{\xi,\zeta}^x$  and $\overline Y_{\xi,\zeta}^x$ (resp. 
definitions \eqref{Yx}--\eqref{bXy} of  
$X_{\xi,\zeta}^y$  and $\overline X_{\xi,\zeta}^y$)
$y\not=y'$ in the first sum while $x\not=x'$ in the second sum
(in both cases $\xi\leq\zeta$ and $\xi>\zeta$).
\er
\br\label{xi=zeta} 
When the two configurations are equal, $\zeta=\xi$, both 
${\overline Y}_{\xi,\xi}^x=\emptyset$ for all $x\in S$
and $X_{\xi,\xi}^y=\emptyset$ for all $y\in S$. The only nonzero 
coupling rates are thus the diagonal terms 
$G _{\xi,\xi}(x,y;x,y)=\Ga_\xi(x,y)$ so that marginals remain equal.
\er
\epr

\subsection{Proof of Theorem \ref{thm:2}}\label{subsec:proof-thm2}
The above increasing Markovian coupling preserves the ordering 
 between marginals when they are ordered  but leaves 
 them otherwise uncoupled. 
In order to deal with 
unordered configurations  and control their discrepancies, 
we show in the next proposition how to build an 
attractive Markov process out of an increasing one. 
\bp\label{propdecdisc}
Suppose that the process defined by 
\eqref{generator} is  monotone  on $\Omega=\{0,1\}^S$. 
Let $\overline{\calL}$ be an
associated increasing process defined on 
$\Omega\times\Omega$ as in Proposition \ref{propcoupling0}, 
with the coupling rates defined in Proposition \ref{propincreasing}. 
 The operator ${\overline \calL}^D$ defined on  all cylinder 
 functions on $\Om\times\Om$ as 
\beq
{\overline \calL}^D f(\xi,\zeta)&=& 
\sum_{x_1,y_1 \in S}  \xi(x_1)  (1-\xi(y_1))  \Ga_\xi(x_1,y_1)
\bigl(f(\xi^{x_1,y_1},\zeta) - f(\xi,\zeta)\bigr)\nonumber\\
&+&\sum_{x_2,y_2 \in S}  \zeta(x_2)  (1-\zeta(y_2))  \Ga_\zeta(x_2,y_2)
\bigl(f(\xi,\zeta^{x_2,y_2}) - f(\xi,\zeta)\bigr)\nonumber\\
&+&\sum_{x_1,y_1\in S}  \sum_{x_2,y_2\in S}
 \xi(x_1)  (1-\xi(y_1))  \zeta(x_2)  (1-\zeta(y_2)) 
 G^D_{\xi,\zeta}(x_1,y_1;x_2,y_2)\nonumber\\
\label{gene2}&&\qquad\qquad
\times \bigl(f(\xi^{x_1,y_1},\zeta^{x_2,y_2}) - f(\xi^{x_1,y_1},\zeta) 
- f(\xi,\zeta^{x_2,y_2}) +f(\xi,\zeta)\bigr)
\label{LD}
\eeq
where for all $(\xi,\zeta)\in \Om\times\Om$, all $(x_1,y_1)\in S^2$ 
and all $(x_2,y_2)\in S^2$, 
\beq\label{GD}
G^D_{\xi,\zeta}(x_1,y_1;x_2,y_2)
&=& \sum_{x,y\in S} (\xi\vee\zeta)(x) (1-(\xi\vee\zeta)(y)) \nonumber\\
&&\qquad\times {\displaystyle \frac{1}{N_{\xi,\zeta}(x,y)}}
G_{\xi,\xi\vee\zeta}(x_1,y_1;x,y) \,G_{\xi\vee\zeta,\zeta}(x,y;x_2,y_2)\\
\nonumber\\
\label{coudis1}
N_{\xi,\zeta}(x,y)&=&
\begin{cases}
\Gamma_{\xi\vee\zeta}(x,y)&
\text{ if } \Gamma_{\xi\vee\zeta}(x,y)> 0\cr
1 &  \hbox{ otherwise. }
\end{cases}
\eeq
is an attractive coupling under which discrepancies do not increase.
\ep
\br\label{rk:GD=G}
When the configurations $\xi,\zeta$ are ordered, $\xi\le\zeta$, 
for all $(x_1,y_1,x_2,y_2)\in S^4$
such that $\xi(x_1)(1-\xi(y_1))\zeta(x_2)(1-\zeta(y_2))\not=0$, 
we have 
\be 
G^D_{\xi,\zeta}(x_1,y_1;x_2,y_2)=G_{\xi,\zeta}(x_1,y_1;x_2,y_2)
\ee
so that  ${\overline \calL}^D f(\xi,\zeta)$ in \eqref{LD} reduces to 
${\overline \calL} f(\xi,\zeta)$ in \eqref{coupling}  when marginals are ordered. 
\er
\subsection{Invariant measures}\label{subsec:inv-meas}
In Proposition \ref{propdecdisc} above, the discrepancies are 
proven to be non increasing, 
but  the characterization of the set of invariant measures, Theorem \ref{prop:inv-meas},
requires a bit more, namely the proof that there is a positive probability 
that any pair of discrepancies 
of opposite sign (that is,  the marginals have opposite occupation numbers,
$\xi(x)>\zeta(x)$, $\xi(y)<\zeta(y)$ for some $x$, $y$ in $S$) 
 disappears in finite time under the coupled process. In contradiction with the case 
of simple exclusion process, this requires additional hypotheses 
on the process. However, 
in close connection to SEP, one may consider processes without 
blocking configurations in the  sense of Definition \ref{def:blocking-conf}.
We then have the following:
\bp\label{decrdiscr}
Consider any exclusion process with generator \eqref{generator} 
such that there are no blocking configurations
and $S$ is fully connected in the sense of Definition \ref{def:blocking-conf}.
Whenever the jump rates are such that all inequalities in 
\eqref{eq:1} and \eqref{eq:2} are strict, 
there exists an increasing coupled process under which 
extremal, translation invariant,  
invariant probability measures are supported on the 
set of coupled configurations
$\{(\xi,\zeta):\,\xi\le\zeta\}\cup\{(\xi,\zeta):\,\xi>\zeta\}$.
\ep
 This result is the crucial step  in the determination of the set
$(\mathcal I\cap\mathcal S)_e$, and in proving Theorem \ref{prop:inv-meas}. 
 This theorem is
 analogous  to \cite[Proposition 3.1]{BGRS2} and to \cite[Theorem 5.13]{gs},
 to which we refer  for a full description of this approach. 
 It has  the same (classical skeleton of) proof, although 
 the transition rates in our case depend on more sites than the departure
and arrival sites of a jump. The key point of the proof is to establish that
for the coupled process, all extremal, translation invariant and  
invariant probability measures are supported on the 
set of coupled configurations
$\{(\xi,\zeta):\,\xi\le\zeta\}\cup\{(\xi,\zeta):\,\xi>\zeta\}$,
 and this is given by 
Proposition \ref{decrdiscr}. 
%
%
%

 In the next Section, we apply our results to various simple but non trivial examples.

\section{Applications} 
\label{sec:Applications}

 In this section we illustrate our results through various examples, 
 and check for them monotonicity conditions of Theorem \ref{thm:1}.
Whenever these conditions are fulfilled, 
we construct the coupled 
generators $\overline{\calL}$ and ${\overline \calL}^D$
by applying Propositions \ref{propcoupling0} and \ref{propdecdisc}. 

In  Subsection \ref{subsec:SEP}, we show that in  the case of simple exclusion,
 our construction reduces to basic coupling.
 In Subsection
\ref{subsec:liggett76}
we consider the exclusion process with speed change
 introduced 
by F. Spitzer in \cite{spitzer} and studied by T.M. Liggett
in \cite{L2}.
In this case, we extend 
the range of previously known attractiveness conditions to necessary and sufficient ones.
Finally, in Subsections \ref{subsec:kstep-et-al} and \ref{subsec:gray-griffeath_and-more} 
we introduce and study  models inspired by traffic flows.

 \subsection{Simple exclusion.}\label{subsec:SEP}
For the simple exclusion process (see Remark \ref{rk:SEP}),
 jump rates are independent on the configuration,
\be
\label{whyH=0} \Ga_\zeta(x,y)-\Ga_\xi(x,y)=0
\ee
for all $\xi$, $\zeta$ in $\Om$ and all $x$, $y$ in $S$.

Attractiveness conditions 
\eqref{eq:1}--\eqref{eq:2} 
reduce to non negativity of jump rates and are thus always satisfied.
We show below that the coupling defined in Proposition 
\ref{propdecdisc} reduces to basic coupling in this case.
 In fact, using simple exclusion rates \eqref{rate:SEP}, 
 the jump rates defined through Formula \eqref{G-rates} 
become,  for all $(x_1,y_1,x_2,y_2)\in S^4$ :
\beq\label{G-diag-SEP}
G _{\xi,\zeta}(x_1,y_1;x_2,y_2)=
\begin{cases}
\delta_{x_1,x_2}\,\delta_{y_1,y_2}\, p(x_1,y_1)& \text{ if $\xi\le\zeta$ or $\xi>\zeta$}\cr
0&\text{ otherwise}
\end{cases}
\eeq
 Therefore, the increasing markovian coupling $\overline{\calL}$ 
 defined through Proposition \ref{propincreasing} coincide with 
 basic coupling on configurations with ordered marginals.  Hence we have
\beq\label{coudis2-SEP}
\overline \varphi_{\xi,\xi\vee\zeta} (x,y) &=& 
 \xi(x)  (1-\xi(y)) p(x,y)\\
\label{coudis3-SEP}
\varphi_{\xi\vee\zeta,\zeta} (x,y) &=&
 \zeta(x)  (1-\zeta(y))  p(x,y)
 \eeq
and
\beq
\label{coudis1-SEP}
N_{\xi,\zeta}(x,y)=
\begin{cases}
p(x,y) & \text{ if } p(x,y)>0\cr
1&\text{ otherwise}
\end{cases}
\eeq
 and using \eqref{G-diag-SEP}, 
\beq\nonumber
G^D_{\xi,\zeta}(x_1,y_1;x_2,y_2)&=&\sum_{x,y} (\xi\vee\zeta)(x) (1-(\xi\vee\zeta)(y))\times\,\nonumber\\
&&\quad\frac{1}{N_{\xi,\zeta}(x,y)}\,
G_{\xi,\xi\vee\zeta}(x_1,y_1;x,y) G_{\xi\vee\zeta,\zeta}(x,y;x_2,y_2)\nonumber\\
&=&\delta_{x_1,x_2}\,\delta_{y_1,y_2}\,(\xi\vee\zeta)(x_1) (1-(\xi\vee\zeta)(y_1))\,  p(x_1,y_1)
\label{GD-diag-SEP}\eeq
%
Finally, the generator of the coupling process defined in Proposition \ref{propdecdisc} reads
\beq
{\overline \calL}^D f(\xi,\zeta)
&=& \sum_{x,y\in S}  p(x,y)
 \xi(x)  (1-\xi(y))  \zeta(x)  (1-\zeta(y)) 
  \bigl(f(\xi^{x,y},\zeta^{x,y}) - f(\xi,\zeta)\bigr)\nonumber\\
&+&\sum_{x,y \in S} p(x,y) \xi(x)  (1-\xi(y)) \bigl(1- \zeta(x)  (1-\zeta(y))\bigr)
\bigl(f(\xi^{x,y},\zeta) - f(\xi,\zeta)\bigr)\nonumber\\
&+&\sum_{x,y \in S} p(x,y) \zeta(x)  (1-\zeta(y))  \bigl(1-  \xi(x)  (1-\xi(y))\bigr)
\bigl(f(\xi,\zeta^{x,y}) - f(\xi,\zeta)\bigr)
\label{LD-SEP}
\eeq
 Hence $\overline{\calL}^D$  identifies to
 the basic coupling generator for SEP. 
This comes from the fact that
non zero coupling rates in \eqref{G-diag-SEP} are diagonal, so that the summation
in formula \eqref{GD} reduces here to a single, diagonal, term.
\subsection{Exclusion processes with speed change}\label{subsec:liggett76}
We consider here a family of models, introduced by F. Spitzer 
in his seminal paper \cite{spitzer},
and later studied by T.M. Liggett in \cite[Part II,Sections 1.1, 4.1]{L2}.
as the product of a configuration dependent velocity $c_\xi(x)$ for the particle at site $x$ 
and a configuration independent jump intensity between sites $x$ and $y$. 
This form is particularly interesting in the original context of a lattice gas. 
The jump rates thus read
\be\label{taux:lig76}
\Gamma_\eta(x,y)= q(x,y) c_\eta(x)\ee
where $q:S\times S\to [0,+\infty)$ satisfies 
for all $x\in S$, $q(x,x)=0$ and
\be\label{q:cond}
\sup_{x\in S}\sum_{y\in S}[q(x,y)+q(y,x)]<\infty
\ee
and $c$ satisfies
\be\label{c:cond}
\sup_{x\in S,\eta\in X}c_\eta(x)<\infty;\qquad
\sup_{x\in S}\sum_{y\in S}[c_{\eta^y}(x)-c_\eta(x)]<\infty
\ee
where 
\be
\eta^y(z)=
\begin{cases}
1- \eta(y) & \text{if $z=y$}\cr
\eta(z) & \text{otherwise}
\end{cases}
\ee
In this context, the monotonicity conditions \eqref{eq:1}--\eqref{eq:2} of Theorem \ref{thm:1} read:

For all $y\in S$ such that $\zeta(y)=0$,
       \begin{equation}
       \label{eq:1L}
 \sum_{x\in S} \xi(x) q(x,y)\, \bigl[c_\xi(x) - c_\zeta(x)\bigr]^+ 
\le
       \sum_{x\in S} \zeta(x) (1 - \xi(x) ) \, q(x,y)\,  c_\zeta(x)  
       \end{equation}
       For all $x\in S$ such that $\xi(x)=1$,
\begin{equation}
       \label{eq:2L}
\bigl(\sum_{y\in S} (1 - \zeta(y) ) \, q(x,y) \bigr) \bigl[c_\zeta(x) - c_\xi(x)\bigr]^+  
\le
      \bigl( \sum_{y\in S} \zeta(y) (1 - \xi(y) ) \,q(x,y)\, \bigr)  c_\xi(x)       
     \end{equation}
T.M. Liggett \cite{L2} constructed an increasing coupling under the more stringent condition 
that the speeds are increasing functions,
\be\label{eq:stringent} 
\text{ For all } \xi\le \zeta, \, \forall x \in S,\, c_\xi(x)\le c_\zeta(x)
\ee
In such a case, equations \eqref{eq:1L} are identically verified, 
while equations \eqref{eq:2L} factorize and read

       For all $x\in S$ such that $\xi(x)=1$,
\begin{equation}\label{sce-dec}
\bigl(\sum_{y\in S} (1 - \zeta(y) ) \, q(x,y) \bigr)\,c_\zeta(x)  
\le
      \bigl( \sum_{y\in S} (1 - \xi(y) ) \,q(x,y)\, \bigr)  c_\xi(x)     
     \end{equation}
which imposes a strong bound on the speed increase, since for all $x\in S$, 
the function $\eta \mapsto \sum_{y\in S} (1 - \eta(y) ) \, q(x,y)$ is a bounded 
decreasing function. For any choice of the jump intensity $q\cdot,\cdot)$, 
one can define a monotone exclusion process with speed change, 
whenever the speed functions are chosen in the form
\be\label{eq:c-eta}
c_\eta(x) =\varphi \bigl(\sum_{y\in S} (1 - \eta(y) ) \, q(x,y)\bigr) 
\text{ for all } x\in S\,\text { and all } \eta\in \Omega
\ee
where $u\mapsto \varphi(u)$ is a decreasing function on $\R^+$ 
such that $\varphi(u)\ge \frac{K}{u}$ for some $K>0$ and  all $u>0$.

Beyond these models, it was not clear whether monotone exclusion processes 
with non increasing speed could exist. We provide below an example with 
decreasing speeds, that is, for all $\xi\le \zeta$, for all $x \in S$, 
$c_\xi(x)\ge c_\zeta(x)$. In such cases, equations \eqref{eq:2L} 
are identically verified, while equations \eqref{eq:1L} now read:

For all $y\in S$ such that $\zeta(y)=0$,
       \begin{equation}
\label{eq:1Ld}
 \sum_{x\in S} \xi(x) q(x,y)\, c_\xi(x)  
\le
       \sum_{x\in S} \zeta(x) \, q(x,y)\,  c_\zeta(x)  
       \end{equation}
so that speed functions have to fulfil a set of coupled inequalities 
indexed by the possible values of $y$.
We take $S=\Z$, fix $L \in \N\setminus\{0\}$ and  set:
\beq\label{eq:c-dec}
\forall (x,y)\in S^2,\- q(x,y) ={\bf 1}_{\{0<y-x\le L\}}\\
\forall \eta\in \Omega, \forall x\in S, \- c_\eta(x) = 2 L - \eta(x)\eta(x+1)
\eeq
Clearly speeds are decreasing functions so that equations \eqref{eq:2L} 
are identically verified.
We now prove that equations \eqref{eq:1Ld} are also fulfilled. 

For any $y\in S$, we have the bound
\be\label{eq:bound}
(2 L-1)\sum_{x=y-L}^{y-1} \eta(x)\le  \sum_{x=y-L}^{y-1} \eta(x) \, 
q(x,y)\,  c_\zeta(x)  \le (2 L+1)\sum_{x=y-L}^{y-1} \eta(x)
\ee
Now for $\xi\le\zeta$, either for all $x\in [y-L, y-1]$, $c_\xi(x) =c_\zeta(x)$ 
and equation \eqref{eq:1Ld} is fulfilled
since $\xi\le\zeta$, or there is $x\in [y-L, y-1]$ such that $\xi(x)=0$ and 
$\zeta(x)=1$, so that, using the previous bounds \eqref{eq:bound},
\be\label{eq:for-1Ld} 
 \sum_{x=y-L}^{y-1} \eta(x) \, q(x,y)\,  
 c_\eta(x) -  \sum_{x=y-L}^{y-1} \eta(x) \, q(x,y)\,  c_\xi(x)  \ge - 2 (L -1) +2 L \ge 2 >0
\ee
and equations \eqref{eq:1Ld} are verified. The exclusion process 
with decreasing speeds defined by \eqref{sce-dec} is monotone.

\subsection{$k$-step exclusion process and related models}\label{subsec:kstep-et-al}
The $k$-step exclusion process was introduced in \cite{guiol}
as an auxiliary model to study the long range exclusion process 
(see also \cite{AG, guiol2}). 
It generalizes the simple exclusion process, we study this model
in dimension 1, when $k=2$, in Subsection \ref{sec:kstep}. We then introduce 
in Subsection \ref{sec:2*step} a 
first variation of the latter model, that we call
$2^*$-step exclusion process. Finally, 
in Subsection \ref{sec:2step-traffic}, we combine both models
to build and analyse a traffic model that we call a range 2 traffic model.  
\subsubsection{The one-dimensional $k$-step exclusion process}\label{sec:kstep}
The state space of the $k$-step exclusion process
is $\{0,1\}^{\Z^d}$. Its jumps follow a translation invariant
probability transition on $\Z^d$. 
In words, if a particle on site $x$ tries to jump, it follows for at most $k$ steps a
random walk $(X^x_n)_{n\geq 0}$ with $X^x_0=x$ until it finds an empty site $y$; 
if all the sites encountered during the $k$ steps are occupied, the particle stays 
on $x$. The  generator of the one-dimensional $k$-step exclusion process is given by 
\begin{equation} \label{eq:kstepgenerator}
\mathfrak{L}_{k}f(\eta )
=\sum_{j=1}^k\sum_{x,y\in\Z} \eta(x)(1-\eta(y)) c_j(x,y,\eta )\left[
f(\eta^{x,y})-f(\eta )\right],
\end{equation}
where $c_j(x,y,\eta)=
{\bf E}^x\left[\prod_{i=1}^{j-1}\eta (X_i),\sigma_y=j\leq
\sigma_x\right] $
and $\sigma_y=\inf
\left\{ n\geq 1:X^x_n=y\right\} $ is the first (non zero) arrival time at site $y$.\\
For the sake of simplicity, we restrict ourselves to the particular case
of the totally asymmetric nearest-neighbor 2-step exclusion
on $S=\Z$, for which we have
\begin{equation} \label{eq:TAnn-2step-rate}
\sum_{j=1}^2 c_j(x,y,\eta )=1_{\{y=x+1\}}+1_{\{y=x+2\}}\eta(x+1 )=:\Ga_\eta(x,y)
\end{equation} 
The totally asymmetric nearest-neighbor 2-step exclusion 
is attractive,   and, as for the simple exclusion process,
 the set $(\mathcal I\cap\mathcal S)_e$ of extremal translation invariant 
 and invariant probability measures for the dynamics 
  consists of a one parameter family $\{\nu_\rho, \rho\in[0,1]\}$
 of Bernoulli product measures, where $\rho$ represents the average density per site,
   see \cite{guiol}. 
This process is a particular case of the range 2 traffic model
studied in subsection \ref{sec:2step-traffic}, hence
its coupling rates are derived as a particular case of Proposition 
\ref{prop:coupl-rates_2step-traffic} below. 
\subsubsection{The one-dimensional totally asymmetric 
$2^*$-step exclusion process}\label{sec:2*step}
On $S=\Z$,  we define
\begin{equation} \label{eq:TAnn-2*step-rate}
\Ga_\eta(x,y)=1_{\{y=x+1\}}+1_{\{y=x+2\}}(1-\eta(x+1 ))
\end{equation} 
 We call \textit{totally asymmetric $2^*$-step exclusion process}
the exclusion process with generator \eqref{generator} for the rate $\Ga_\eta(x,y)$
given in \eqref{eq:TAnn-2*step-rate}.  
The totally asymmetric nearest-neighbor 2*-step exclusion is
 a particular case of the range 2 traffic model
studied in Subsection \ref{sec:2step-traffic} below, 
hence its attractiveness follows from
Proposition \ref{lem:2step-traffic-attra}, and its coupling rates
 are derived as a particular case of Proposition 
\ref{prop:coupl-rates_2step-traffic}. \\ 
This model is also a particular case of a more general 
$2^*$-step exclusion process of transition rate  given by
\begin{equation} \label{eq:2*step-rate}
\Ga_\eta(x,y)=p(x,y)+\sum_{z\in\Z}p(x,z)p(z,y)(1-\eta(z))
\end{equation} 
for a translation invariant
transition probability $p(.,.)$.
\begin{proposition}\label{prop:inv-2*step}
The Bernoulli product measures $\{\nu_\rho, \rho\in[0,1]\}$ are 
invariant for the 2*-step exclusion process of transition rate 
$\Ga_\eta(x,y)$ given in \eqref{eq:2*step-rate}.
\end{proposition}
\bpr
We proceed as in the proof of \cite[Theorem VIII.2.1]{L}, by checking that
$\int L f_A d\nu_\rho=0$, where $A$ is a finite set of sites
and $f_A$ is the cylinder function defined by
\begin{equation}\label{eq:fA}
f_A(\eta)=\prod_{x\in A} \eta(x)
\end{equation}
We have, denoting by $\calL_{\rm SEP}$ the generator of the
simple exclusion process and by $\calL_{\rm 2*s}$ the second part of
the generator of the 
2*-step exclusion process,
\begin{eqnarray*}\nonumber
\calL f_A(\eta)&=&\calL_{\rm SEP} f_A(\eta)+\calL_{\rm 2*s} f_A(\eta)\\
\calL_{\rm 2*s} f_A(\eta)&=&
\sum_{x,y\in S,x\not=y}\sum_{z\in S,z\not=x,y} p(x,z)p(z,y)
\eta(x)(1-\eta(y))(1-\eta(z)) \bigl[f_A(\eta^{x,y})-f_A(\eta)\bigr]\label{eq:gen-fA}
\end{eqnarray*}
Since
\begin{equation*}
\int f_A(\eta) \eta(x) (1-\eta(z)) (1-\eta(y)) d\nu_\rho(\eta)=
\begin{cases}
 0, & \text{ if } y\in A \text{ or } z\in A\\
(1-\rho)^2 \rho^{|A\cup\{x\}|},
 & \text{ if } y\notin A, z\notin A
\end{cases}
\end{equation*}
and
\begin{equation*}
\int f_A(\eta^{x,y}) \eta(x) (1-\eta(z)) (1-\eta(y)) d\nu_\rho(\eta)=
\begin{cases}
 0, & \text{ if } x\in A \text{ or } z\in A\\
(1-\rho)^2  \rho^{|A\cup\{x\}\setminus\{y\}|}, 
 & \text{ if } x\notin A, z\notin A
\end{cases}
\end{equation*}
we have
\begin{eqnarray*}
\int \calL_{\rm 2*s}f_A(\eta)d\nu_\rho(\eta)&=&
\sum_{x,y:x\not=y,x\notin A}\sum_{z:z\not=x,y,z\notin A}
p(x,z)p(z,y)(1-\rho)^2 \rho^{|A\cup\{x\}\setminus\{y\}|}\\
&&-\sum_{x,y:x\not=y,y\notin A}\sum_{z:z\not=x,y,z\notin A}
p(x,z)p(z,y)(1-\rho)^2 \rho^{|A\cup\{x\}|}
\end{eqnarray*}
Taking $x\notin A,y\notin A$ in the first sums of the two terms 
on the right hand side gives 0, hence
we are left with $y\in A$ for the first term, and
$x\in A$ for the second term. Exchanging the indexes $x$ and $y$
in the second term gives
\begin{eqnarray*}
\int \calL_{\rm 2*s}f_A(\eta)d\nu_\rho(\eta)&=&
(1-\rho)^2 \rho^{|A|}
\sum_{x,y:x\not=y,x\notin A,y\in A}\sum_{z:z\not=x,y,z\notin A}
\bigl[p(x,z)p(z,y)-
p(y,z)p(z,x)\bigr]\\
&=& 0
\end{eqnarray*}
because $A$ is finite and $p(.,.)$ is bi-stochastic.
\epr
 \subsubsection{A range 2 traffic model}\label{sec:2step-traffic}
On $S=\Z$, for $\alpha,\beta\in[0,1]$, we define
\begin{equation} \label{eq:Traffic-2step-rate}
\Ga_\eta(x,y)=1_{\{y=x+1\}}+1_{\{y=x+2\}}[\alpha\eta(x+1 )+\beta (1-\eta(x+1 ))]
\end{equation} 
We call \textit{range 2 traffic model}
the exclusion process with generator \eqref{generator} for the rate $\Ga_\eta(x,y)$
given in \eqref{eq:Traffic-2step-rate}. 
This rate is a convex combination of the respective rates for one-dimensional 
totally asymmetric simple exclusion, 2-step exclusion and $2^*$-step exclusion. 
The traffic interpretation is that a car can either go one step ahead, 
or 2 steps ahead by overtaking another car or by accelerating.
\begin{proposition}\label{lem:2step-traffic-attra}
The range 2 traffic model is attractive if and only if $|\beta-\alpha|\leq 1$.
The case $\beta=\alpha=0$ corresponds to simple exclusion, the case $\beta=0$
to 2-step exclusion, and the case $\alpha=0$ to $2^*$-step exclusion.
\end{proposition}
\bpr
We have to check inequalities \eqref{eq:1}--\eqref{eq:2}. Let $(\xi,\zeta)\in \Omega^2$
be such that $\xi\le\zeta$.

We begin with \eqref{eq:1}. Let $y\in \Z$ be such that $\zeta(y)=0$, hence 
 $\xi(y)=0$. Then we write
 \begin{eqnarray}\nonumber
 &\sum_{x\in\Z}& \zeta(x) (1 - \xi(x) ) \Ga_\zeta(x,y)\\
\label{eq:lhs-eq1} &&=
  \zeta(x-1)(1 - \xi(x-1)) + \zeta(x-2) (1 - \xi(x-2) )(\alpha\zeta(x-1) + \beta (1 - \zeta(x-1))\\\nonumber
    &\sum_{x\in\Z}& \xi(x) \bigl[\Ga_\xi(x,y) - \Ga_\zeta(x,y)\bigr]^+ \\\label{eq:rhs-eq1} &&=
\xi(x-2)\bigl[ (\alpha\xi(x-1) + \beta (1 - \xi(x-1))- (\alpha\zeta(x-1)
 + \beta (1 - \zeta(x-1))   \bigr]^+    
     \end{eqnarray}
First, if $\xi(x-1)=1$ then $\zeta(x-1)=1$, hence \eqref{eq:rhs-eq1} 
is null; secondly, if $\zeta(x-1)=0$ then $\xi(x-1)=0$, hence \eqref{eq:rhs-eq1} 
is null; in both cases, \eqref{eq:1} is satisfied. Finally, if $\xi(x-1)=0$ and 
$\zeta(x-1)=1$, then \eqref{eq:lhs-eq1} is equal to 
$1+\alpha\zeta(x-2) (1 - \xi(x-2) )$ while \eqref{eq:rhs-eq1} 
is equal to $\xi(x-2) ( \beta - \alpha)^+$: either $\xi(x-2)=0$ and 
\eqref{eq:1} is satisfied, or $\xi(x-2)=1$ and $( \beta - \alpha)^+\le 1$ 
is required for \eqref{eq:1} to be satisfied.\\

We now check \eqref{eq:2}. Let $x\in \Z$ be such that $\xi(x)=1$, hence 
 $\zeta(x)=1$. Then we write
\begin{eqnarray}\nonumber      
      & \sum_{y\in\Z}& \zeta(y) (1 - \xi(y) ) \Ga_\xi(x,y) \\  \label{eq:lhs-eq2}
 &&= \zeta(x+1) (1 - \xi(x+1) ) + \zeta(x+2) (1 - \xi(x+2) )(\alpha\xi(x+1) 
 + \beta (1 - \xi(x+1))
  \\\nonumber
  &\sum_{y\in\Z}& (1 - \zeta(y) )\bigl [\Ga_\zeta(x,y) - \Ga_\xi(x,y)\bigr]^+  \\\label{eq:rhs-eq2}
 &&= (1 - \zeta(x+2) )\bigl [(\alpha\zeta(x+1) 
 + \beta (1 - \zeta(x+1))-(\alpha\xi(x+1) + \beta (1 - \xi(x+1))\bigr]^+ 
        \end{eqnarray}
 First, if $\xi(x+1)=1$ then $\zeta(x+1)=1$, hence \eqref{eq:rhs-eq2} 
 is null; secondly, if $\zeta(x+1)=0$ then $\xi(x+1)=0$, hence 
 \eqref{eq:rhs-eq2} is null; in both cases, \eqref{eq:2} is satisfied. 
 Finally, if $\xi(x+1)=0$ and $\zeta(x+1)=1$, then \eqref{eq:lhs-eq2} 
 is equal to $1+\beta\zeta(x+2) (1 - \xi(x+2) )$ while \eqref{eq:rhs-eq2} 
is equal to $(1-\zeta(x+2)) ( \alpha - \beta)^+$: either $\zeta(x+2)=1$ 
and \eqref{eq:2} is satisfied, or $\zeta(x+2)=0$ and 
$( \alpha - \beta)^+\le 1$ is required for \eqref{eq:2} to be satisfied.
\epr

\noindent \textbf{Invariant measures.}
Because it is the case for simple exclusion, 2-step exclusion 
and 2*-step exclusion processes (see Proposition \ref{prop:inv-2*step})
the Bernoulli product measures $\{\nu_\rho, \rho\in[0,1]\}$ are 
invariant for the range 2 traffic model. This model has no blocking configurations
if $\alpha,\beta$ are positive, in which case the Bernoulli product measures are
the extremal translation invariant and invariant probability measures for the dynamics, 
by Theorem \ref{prop:inv-meas}. \\

\noindent \textbf{Coupling rates.}
 To compute the coupling rates for dynamics on $S=\Z$, we use
the following formulas, equivalent to \eqref{G-rates}.
For all $(x_1,y_1)\in S^2$, 
\be\label{G-diag}
G _{\xi,\zeta}(x_1,y_1;x_1,y_1) =\Ga_\xi(x_1,y_1)\wedge\Ga_\zeta(x_1,y_1) 
\ee
For all $(x_1,y_1,x_2,y_2)\in S^4$ such that $(x_1,y_1)\not=(x_2,y_2)$,
\beq\nonumber
G _{\xi,\zeta}(x_1,y_1;x_2,y_2) &=&
\begin{cases}
\delta_{x_1,x_2} \bigl[H^i_{\xi,\zeta}(x_1;y_1,y_2)\bigr]^+ 
+ \delta_{y_1,y_2} \bigl[H^f_{\xi,\zeta}(x_1,x_2;y_1)\bigr]^+
& \hbox{ if } \xi\le \zeta \cr
\delta_{x_1,x_2} \bigl[H^i_{\zeta,\xi}(x_1;y_2,y_1)\bigr]^+ 
+ \delta_{y_1,y_2} \bigl[H^f_{\zeta,\xi}(x_2,x_1;y_1)\bigr]^+
& \hbox{ if } \xi > \zeta \cr
0& \hbox{ otherwise. }
\end{cases}\\\label{G-hors-diag}
\eeq
with,  for all $(x,y,z)\in S^3$, 
\beq\label{Hi}
H^i_{\xi,\zeta}(x;y,z) &=&  \bigl(\sum_{y'\le y} (1-\xi(y'))\zeta(y')\Ga_\xi(x,y')\bigr)\wedge 
 \bigl(\sum_{z'\le z} (1-\zeta(z'))\bigl[\Ga_\zeta(x,z')-\Ga_\xi(x,z')\bigr]^+\bigr)\nonumber\\
 &&\qquad-   \bigl(\sum_{y'< y} (1-\xi(y'))\zeta(y')\Ga_\xi(x,y')\bigr)\vee 
  \bigl(\sum_{z'< z} (1-\zeta(z'))\bigl[\Ga_\zeta(x,z')-\Ga_\xi(x,z')\bigr]^+\bigr)\\
  \label{Hf}
 H^f_{\xi,\zeta}(x,y;z) &=&\bigl(\sum_{x'\le x} \xi(x')\bigl[\Ga_\xi(x',z)-\Ga_\zeta(x',z)\bigr]^+\bigr)
 \wedge 
 \bigl(\sum_{y'\le y} \zeta(y')(1-\xi(y'))\Ga_\zeta(y',z)\bigr)\nonumber\\
 &&\qquad \quad-\bigl(\sum_{x'< x} \xi(x')\bigl[\Ga_\xi(x',z)-\Ga_\zeta(x',z)\bigr]^+\bigr)
 \vee
 \bigl(\sum_{y'< y} \zeta(y')(1-\xi(y'))\Ga_\zeta(y',z)\bigr)
  \eeq 
Therefore, for the range 2 traffic model, 
applying Propositions \ref{propcoupling0}, \ref{propincreasing}, 
 and formulas \eqref{G-diag}--\eqref{Hf}, 
we obtain first the following formulas for the coupling rates
$G _{\xi,\zeta}(x_1,y_1;x_2,y_2)$, taking into account that 
in formula \eqref{coupling}, they are multiplied by the prefactor 
 $\xi(x_1)(1-\xi(y_1))\zeta(x_2)(1-\zeta(y_2))$,
 so that $\xi(x_1)=1-\xi(y_1)=\zeta(x_2)=1-\zeta(y_2)=1$:
\begin{eqnarray}\nonumber
 G_{\xi,\zeta}(x,x+1;x,x+1) &=&1\\\nonumber
 G_{\xi,\zeta}(x,x+2;x,x+2) &=&[\alpha\xi(x+1 )+\beta (1-\xi(x+1 ))]\wedge[\alpha\zeta(x+1 )+\beta (1-\zeta(x+1 ))]\\
 G_{\xi,\zeta}(x,x+1;x,x+2) &=& 
[\alpha-\beta]^+\zeta(x+1)
 \quad\hbox{ when  }\xi\le\zeta\nonumber\\ 
 &=& 0
 \quad\hbox{ when  }\xi>\zeta\nonumber\\
 G_{\xi,\zeta}(x,x+2;x+1,x+2) &=& [\beta-\alpha]^+(1-\xi(x+1))\quad\hbox{ when  }\xi\ge\zeta\nonumber\\ 
 &=& 0
 \quad\hbox{ when  }\xi>\zeta\nonumber\\
 G_{\xi,\zeta}(x+1,x+2;x,x+2)&=& 0\quad\hbox{ when  }\xi\ge\zeta\nonumber\\ 
 &=& (\beta-\alpha)^+(1-\zeta(x+1))
 \quad\hbox{ when  }\xi>\zeta\nonumber\\
  G_{\xi,\zeta}(x,x+2;x,x+1) &=& 0\quad\hbox{ when  }\xi\ge\zeta\nonumber\\ 
 &=& 0
 \quad\hbox{ when  }\xi>\zeta\nonumber\\
 G_{\xi,\zeta}(x+1,x+2;x,x+2)&=& 0\quad\hbox{ when  }\xi\ge\zeta\nonumber\\ 
 &=& [\alpha-\beta]^+\xi(x+1)
 \quad\hbox{ when  }\xi>\zeta\nonumber
\end{eqnarray}
Some more computations to get the formulas in Proposition 
\ref{propdecdisc} yield:
\begin{proposition}\label{prop:coupl-rates_2step-traffic}
The  coupled generator of the range 2 traffic model writes
\begin{equation}\label{LD_2step-traffic}
{\overline \calL}^D f(\xi,\zeta)={\overline \calL}^D_1 f(\xi,\zeta)+
{\overline \calL}^D_2 f(\xi,\zeta)+{\overline \calL}^D_{3,1} f(\xi,\zeta)
+{\overline \calL}^D_{3,2} f(\xi,\zeta)
\end{equation}
where ${\overline \calL}^D_1$ deals with coupled jumps with the same departure 
and arrival sites, ${\overline \calL}^D_2$ with coupled jumps with a 
different site either for departure or for arrival, and ${\overline \calL}^D_{3,1}$, 
${\overline \calL}^D_{3,2}$ deal with uncoupled jumps. They are given by
\begin{eqnarray}
{\overline \calL}^D_1 f(\xi,\zeta)&=& 
\sum_{x \in S} \xi(x)  (1-\xi(x+1))  \zeta(x)  (1-\zeta(x+1))
\times \bigl(f(\xi^{x,x+1},\zeta^{x,x+1}) - f(\xi,\zeta)\bigr)\nonumber\\
&+&\sum_{x \in S} \xi(x)\zeta(x)(1-(\xi\vee\zeta)(x+2)) \nonumber\\
&&\quad\times\bigl[\alpha\xi(x+1)\zeta(x+1)+\beta(1-(\xi\vee\zeta)(x+1))\bigr.\nonumber\\
&&\qquad\bigl. +(\alpha\wedge\beta)\{\xi(x+1)(1-\zeta(x+1)) +\zeta(x+1)(1-\xi(x+1))\}
\bigr]\nonumber\\
&&\quad\times \bigl(f(\xi^{x,x+2},\zeta^{x,x+2}) - f(\xi,\zeta)\bigr)
\label{LD1_2step-traffic}
\end{eqnarray}
\begin{eqnarray}
{\overline \calL}^D_2 f(\xi,\zeta)&=& 
\sum_{x \in S} \xi(x)(1-\xi(x+1))\zeta(x+1)(1-(\xi\vee\zeta)(x+2))(\beta-\alpha)^+ \nonumber\\
&&\quad
\times \bigl(f(\xi^{x,x+2},\zeta^{x+1,x+2}) - f(\xi,\zeta)\bigr)\nonumber\\
&+&\sum_{x \in S} \xi(x)(1-\zeta(x))\zeta(x-1)(1-(\xi\vee\zeta)(x+1))
(\beta-\alpha)^+\nonumber\\
&&\quad
\times \bigl(f(\xi^{x,x+1},\zeta^{x-1,x+1}) - f(\xi,\zeta)\bigr)\nonumber\\
&+&\sum_{x \in S} \xi(x)\zeta(x)\xi(x+1)(1-\zeta(x+1))
(1-(\xi\vee\zeta)(x+2))(\alpha-\beta)^+ \nonumber\\
&&\quad
\times
\bigl(f(\xi^{x,x+2},\zeta^{x,x+1}) - f(\xi,\zeta)\bigr)\nonumber\\
&+&\sum_{x \in S} \xi(x)\zeta(x)(1-\xi(x+1))\zeta(x+1)
(1-(\xi\vee\zeta)(x+2))(\alpha-\beta)^+\nonumber\\
&&\quad\times \bigl(f(\xi^{x,x+1},\zeta^{x,x+2}) - f(\xi,\zeta)\bigr)
\label{LD2_2step-traffic}
\end{eqnarray}
\begin{eqnarray}
{\overline \calL}^D_{3,1} f(\xi,\zeta)&=& 
\sum_{x \in S} \xi(x) (1-\xi(x+1)) \bigl[ 1- \zeta(x)  (1-\zeta(x+1))\nonumber\\
&&\quad
-(1-\zeta(x))\zeta(x-1)(1-\zeta(x+1))(\beta-\alpha)^+\bigr.\nonumber\\
&&\qquad\bigl.
-\zeta(x)\zeta(x+1)(1-(\xi\vee\zeta)(x+2))(\alpha-\beta)^+
\bigr]\nonumber\\
&&\quad
\times
\bigl(f(\xi^{x,x+1},\zeta) - f(\xi,\zeta)\bigr)\nonumber\\
&+&\sum_{x \in S} \xi(x) (1-\xi(x+2)) 
\bigl[ \alpha\xi(x+1)+\beta(1-\xi(x+1))\bigr.\nonumber\\
&&\quad\bigl.   - \zeta(x) (1-\zeta(x+2))\bigl\{
\alpha\xi(x+1)\zeta(x+1)+\beta(1-(\xi\vee\zeta)(x+1))\bigr.\nonumber\\
&&\qquad\qquad\bigl. +(\alpha\wedge\beta)\{\xi(x+1)(1-\zeta(x+1)) +\zeta(x+1)(1-\xi(x+1))\}
\bigr\}\nonumber\\
&&\quad\bigl. -\zeta(x)\xi(x+1)(1-\zeta(x+1))(1-\zeta(x+2))(\alpha-\beta)^+
\bigr]\nonumber\\
&&\quad
\times
\bigl(f(\xi^{x,x+2},\zeta) - f(\xi,\zeta)\bigr)
\label{LD31_2step-traffic}
\end{eqnarray}
\begin{eqnarray}
{\overline \calL}^D_{3,2} f(\xi,\zeta)&=& 
\sum_{x \in S} \zeta(x) (1-\zeta(x+1)) 
\bigl[1- \xi(x)(1-\xi(x+1))\bigr.\nonumber\\
&&\qquad\bigl. 
 - \xi(x)\xi(x+1) (1-(\xi\vee\zeta)(x+2))(\alpha-\beta)^+\bigr.\nonumber\\
&&\qquad\bigl. -(1-\xi(x))\xi(x-1)(1-\xi(x+1))(\beta-\alpha)^+
\bigr]\nonumber\\
&&\quad
\times
\bigl(f(\xi,\zeta^{x,x+1}) - f(\xi,\zeta)\bigr)\nonumber\\
&+&\sum_{x \in S} \zeta(x) (1-\zeta(x+2)) 
\bigl[ \alpha\zeta(x+1)+\beta(1-\zeta(x+1))\bigr.\nonumber\\
&&\qquad\bigl.   - \xi(x) (1-\xi(x+2))\bigl\{
\alpha\xi(x+1)\zeta(x+1)+\beta(1-(\xi\vee\zeta)(x+1))
\bigr.\nonumber\\
&&\qquad\qquad\bigl. +(\alpha\wedge\beta)\{\zeta(x+1)(1-\xi(x+1)) +\xi(x+1)(1-\zeta(x+1))
\bigr\}\nonumber\\
&&\qquad\bigl.
-\xi(x)(1-\xi(x+1))\zeta(x+1)(1-\xi(x+2))(\alpha-\beta)^+\bigr]\nonumber\\
&&\quad\times
\bigl(f(\xi,\zeta^{x,x+2}) - f(\xi,\zeta)\bigr)\label{LD32_2step-traffic}
\end{eqnarray}
\end{proposition} 
\begin{remark}\label{rk:sep-2step-2step*_coupl}
Taking $\alpha=\beta=0$ gives the basic coupling  generator
for TASEP, while taking  $\beta=0$ gives a coupled generator 
for $2^*$-step exclusion, and taking $\alpha=0$ gives a coupled generator 
for 2-step exclusion. The latter is different from the one used in \cite{guiol},
where coupled jumps had the same departure site, 
but there were no coupled jumps with the same arrival site but with different departure sites.
\end{remark}  
\subsection{From a non-attractive traffic model to an attractive dynamics}\label{subsec:gray-griffeath_and-more}
We begin with an exclusion process with the
transition rates introduced in \cite{GG}
in the context of a cellular automaton dynamics. There, $S=\Z$, and the transitions
are nearest neighbor and totally asymmetric. For all $x\in S$,
$\eta\in X$ such that $\eta(x)=1$ and $\eta(x+1)=0$
\begin{eqnarray}
\label{taux:gg}
\Gamma_\eta(x,x+1)=
  \begin{cases}
    \alpha & \text{ if } \eta(x-1)=1,\, \eta(x+2)=0, \text{ [accelerating] }\\
    \beta & \text{ if } \eta(x+2)=1,\, \eta(x-1)=0, \text{ [braking] }\\
    \gamma & \text{ if } \eta(x-1)=\eta(x+2)=1, \text{ [congested] }\\
    \delta & \text{ if }\eta(x-1)=\eta(x+2)=0, \text{ [driving] }.
  \end{cases}
\end{eqnarray}
 where the parameters  $\alpha,\beta, \gamma, \delta$ are positive.
 This model is not attractive, unless it reduces to simple exclusion, that is
$\alpha=\beta=\gamma=\delta$.   Indeed, for any other choice, conditions \eqref{eq:1}--\eqref{eq:2} 
from Theorem \ref{thm:1} are not  satisfied.
Here, it is possible to turn the dynamics into an attractive one, just by considering a symmetrized version,
in which the non zero rates are the previous, rightwards, ones, \eqref{taux:gg}, together with the following symmetric, leftwards rates:
\begin{eqnarray}
\label{taux:symmetric-for-gg}
\Gamma_\eta(x+1,x)=
  \begin{cases}
    \alpha & \text{ if } \eta(x+2)=1,\, \eta(x-1)=0, \\
    \beta & \text{ if } \eta(x-1)=1,\, \eta(x+2)=0, \\
    \gamma & \text{ if } \eta(x+2)=\eta(x-1)=1, \\
    \delta & \text{ if }\eta(x+2)=\eta(x-1)=0.
  \end{cases}
\end{eqnarray}
Then applying conditions in Theorem \ref{thm:1}  leads to the following result.
\bp\label{prop:attra-gg-extended}
The symmetrized dynamics with rates \eqref{taux:gg}--\eqref{taux:symmetric-for-gg}
is attractive if and only if $\alpha,\beta, \gamma, \delta$ satisfy the following
conditions
\be\label{cond-for-gg-extended}
\beta\leq\gamma\wedge\delta\leq\gamma\vee\delta\leq\alpha,\qquad
\alpha\leq\beta+\gamma\wedge\delta,\qquad\delta\leq 2\beta
\ee
\ep
\noindent
Note that the facilitated exclusion process (\cite{AGLS,BBCS,BESS})
has rates \eqref{taux:gg} with $\alpha=\gamma=1,\, \beta=\delta=0$.
Hence it is not attractive, and its symmetrized version 
(with the corresponding rates in \eqref{taux:symmetric-for-gg}) 
is not attractive either. Indeed the study of this model required
other tools. 

\noindent \textbf{Invariant measures.}
The symmetrized dynamics has no blocking configurations, since the parameters
$\alpha,\beta, \gamma, \delta$ are positive. We can
thus apply Theorem \ref{prop:inv-meas}.

 \noindent \textbf{Coupling rates.}
Applying Propositions \ref{propcoupling0}, \ref{propincreasing}, 
with formulas \eqref{G-diag}--\eqref{Hf}, 
we obtain first the following formulas for the coupling rates
$G _{\xi,\zeta}(x_1,y_1;x_2,y_2)$, taking into account that  
they are multiplied by the prefactor 
 $\xi(x_1)(1-\xi(y_1))\zeta(x_2)(1-\zeta(y_2))$,
 so that $\xi(x_1)=1-\xi(y_1)=\zeta(x_2)=1-\zeta(y_2)=1$.
 Note that since the rates \eqref{taux:gg} and 
 \eqref{taux:symmetric-for-gg} are symmetric, it is enough to compute
 the coupling rates in the positive direction to get the ones 
 in the negative direction by symmetry. 
 To simplify the computations, we assume that $\gamma\leq\delta$.
\begin{eqnarray}\nonumber
 G_{\xi,\zeta}(x,x+1;x,x+1) &=&\Gamma_\xi(x,x+1)\wedge\Gamma_\zeta(x,x+1)\\\nonumber
 G_{\xi,\zeta}(x,x-1;x,x-1) &=&\Gamma_\xi(x,x-1)\wedge\Gamma_\zeta(x,x-1)
 \end{eqnarray}
When $\xi\leq\zeta$, we have
\begin{eqnarray}\nonumber
 G_{\xi,\zeta}(x,x+1;x,x-1) &=&\zeta(x+1)
 \bigl[(1-\zeta(x-2))(\alpha-\delta) + \xi(x-2)(\gamma -\beta)\bigr]\\\nonumber
 G_{\xi,\zeta}(x,x-1;x,x+1) &=&\zeta(x-1)
 \bigl[(1-\zeta(x+2))(\alpha-\delta) + \xi(x+2)(\gamma -\beta)\bigr]\\\nonumber
G_{\xi,\zeta}(x,x+1;x+2,x+1) &=& (1-\xi(x+2))\bigl[(1-\zeta(x-1))(\delta-\beta)\\\nonumber
 &&\quad
 +(1- \xi(x-1))\zeta(x-1)(\delta-\gamma)+\xi(x-1)(\alpha-\gamma)\bigr]\\\nonumber
 G_{\xi,\zeta}(x,x-1;x-2,x-1) &=& (1-\xi(x-2))\bigl[(1-\zeta(x+1))(\delta-\beta)\\\nonumber
 &&\quad
 +(1- \xi(x+1))\zeta(x+1)(\delta-\gamma)+\xi(x+1)(\alpha-\gamma)\bigr]
  \end{eqnarray}
When $\xi>\zeta$, we have
\begin{eqnarray}\nonumber
 G_{\xi,\zeta}(x,x+1;x,x-1) &=&\xi(x-1)
\bigl[(1- \xi(x+2))(\alpha-\delta)+ \zeta(x+2)(\gamma-\beta)\bigr]\\\nonumber
 G_{\xi,\zeta}(x,x-1;x,x+1) &=&\xi(x+1)
\bigl [(1- \xi(x-2))(\alpha-\delta)+ \zeta(x-2)(\gamma-\beta)\bigr]\\\nonumber
  G_{\xi,\zeta}(x,x+1;x+2,x+1) &=& (1-\zeta(x))\bigl[(1-\xi(x+3))(\delta-\beta)\\\nonumber
 &&\quad
 +\zeta(x+3)(\alpha-\gamma)+\xi(x+3)(1-\zeta(x+3))(\delta-\gamma)\bigr]\\\nonumber
  G_{\xi,\zeta}(x,x-1;x-2,x-1) &=& (1-\zeta(x))\bigl[(1-\xi(x-3))(\delta-\beta)\\\nonumber
 &&\quad
 +\zeta(x-3)(\alpha-\gamma)+\xi(x-3)(1-\zeta(x-3))(\delta-\gamma)\bigr]
\end{eqnarray}
Finally, applying Proposition \ref{propdecdisc}
with formulas \eqref{GD}--\eqref{coudis1}, 
we obtain the following formulas for the coupling rates
$G^D_{\xi,\zeta}(x_1,y_1;x_2,y_2)$, taking into account that  
they are multiplied by the prefactor 
 $\xi(x_1)(1-\xi(y_1))\zeta(x_2)(1-\zeta(y_2))$,
 so that $\xi(x_1)=1-\xi(y_1)=\zeta(x_2)=1-\zeta(y_2)=1$.
Again, since the rates \eqref{taux:gg} and 
 \eqref{taux:symmetric-for-gg} are symmetric, it is enough to compute
 the coupling rates in the positive direction.
 \begin{eqnarray}\nonumber
 G^D_{\xi,\zeta}(x,x+1;x,x+1) &=&\Gamma_\xi(x,x+1)\wedge\Gamma_\zeta(x,x+1)\\\nonumber
 G^D_{\xi,\zeta}(x,x-1;x,x-1) &=&\Gamma_\xi(x,x-1)\wedge\Gamma_\zeta(x,x-1)\\\nonumber
 G^D_{\xi,\zeta}(x,x+1;x,x-1) &=&(1-\zeta(x+1))\xi(x-1)\times\\\nonumber
&&\qquad 
 \bigl[(1-(\xi\vee\zeta)(x+2))(\alpha-\delta) + \zeta(x+2)(\gamma -\beta)\bigr]\\\nonumber
  &&
 +(1-\xi(x-1))\zeta(x+1)\times\\\nonumber
&&\qquad 
 \bigl[(1-(\xi\vee\zeta)(x-1))(\alpha-\delta) + \xi(x-2)(\gamma -\beta)\bigr]\\\nonumber
 G^D_{\xi,\zeta}(x,x-1;x,x+1) &=&(1-\zeta(x-1))\xi(x+1)\times\\\nonumber
&&\qquad 
 \bigl[(1-(\xi\vee\zeta)(x-2))(\alpha-\delta) + \zeta(x-2)(\gamma -\beta)\bigr]\\\nonumber
&& 
 +(1-\xi(x+1))\zeta(x-1)\times\\\nonumber
&&\quad 
 \bigl[(1-(\xi\vee\zeta)(x+2))(\alpha-\delta) + (\xi\vee\zeta)(x+2)(\gamma -\beta)\bigr]\\\nonumber
 G^D_{\xi,\zeta}(x,x+1;x+2,x+1) &=& (1-\zeta(x))
 \bigl[(1-(\xi\vee\zeta)(x+3))(\delta-\beta)\\\nonumber
 &&\quad
 +\zeta(x+3)(\alpha-\gamma)+(1- \zeta(x+3))\xi(x+3)(\delta-\gamma)\bigr]\times\\\nonumber
&&\qquad \bigl[\xi(x+2)\zeta(x-1)(1- \xi(x-1))\left(\frac{\beta}{\gamma}-1\right)+1\bigr]\\\nonumber
 &&
 +(1-\xi(x+2))
 \bigl[(1-(\xi\vee\zeta)(x-1))(\delta-\beta)\\\nonumber
 &&\quad
 +\xi(x-1)(\alpha-\gamma)+(1- \xi(x-1))\zeta(x-1)(\delta-\gamma)\bigr]\times\\\nonumber
&&\qquad \bigl[\xi(x+3)\zeta(x)(1- \zeta(x+3))\left(\frac{\beta}{\gamma}-1\right)+1\bigr]\\\nonumber
 G^D_{\xi,\zeta}(x,x-1;x-2,x-1) &=& (1-\zeta(x))\bigl[(1-(\xi\vee\zeta)(x-3))(\delta-\beta)\\\nonumber
&&\quad
+(1- \zeta(x-3))\xi(x-3)(\delta-\gamma)+\zeta(x-3)(\alpha-\gamma)\bigr]\times\\\nonumber
&&\qquad \bigl[\xi(x-2)\zeta(x+1)(1- \xi(x+1))\left(\frac{\beta}{\gamma}-1\right)+1\bigr]\\\nonumber
 &&
 +(1-\xi(x-2))
 \bigl[(1-(\xi\vee\zeta)(x+1))(\delta-\beta)\\\nonumber
 &&\quad
 +\xi(x+1)(\alpha-\gamma)+(1- \xi(x+1))\zeta(x+1)(\delta-\gamma)\bigr]\times\\\nonumber
&&\qquad \bigl[\xi(x-3)\zeta(x)(1- \zeta(x-3))\left(\frac{\beta}{\gamma}-1\right)+1\bigr]\\\nonumber
G^D_{\xi,\zeta}(x,x+1;x+2,x+3) &=& (1-\zeta(x+1))(1- \xi(x+2))\xi(x+3)\zeta(x)\frac{(\gamma-\beta)}{\gamma}\times\\\nonumber
&&\quad \bigl[(1-(\xi\vee\zeta)(x-1))(\delta-\beta)\\\nonumber
 &&\quad
 +(1- \xi(x-1))\zeta(x-1)(\delta-\gamma)+\xi(x-1)(\alpha-\gamma)\bigr]\\\nonumber
 G^D_{\xi,\zeta}(x,x-1;x-2,x-3) &=& (1-\zeta(x-1))(1- \xi(x-2))\xi(x_3)\zeta(x)\frac{(\gamma-\beta)}{\gamma}\times\\\nonumber
&&\quad \bigl[(1-(\xi\vee\zeta)(x+1))(\delta-\beta)\\\nonumber
 &&\quad
 +(1- \xi(x+1))\zeta(x+1)(\delta-\gamma)+\xi(x+1)(\alpha-\gamma)\bigr]
  \end{eqnarray}
\section{Technical proofs}\label{sec:proofs}
\bpr[Proposition \ref{propcoupling0}]
We rewrite the generator  $\overline\calL$ (Equation \eqref{coupling}) as,
\beq\nonumber
\overline{\calL}f(\xi,\zeta)&=&
\sum_{x_1,y_1 \in S} \xi(x_1)( 1-\xi(y_1))
\left(\Ga_\xi(x_1,y_1)
-  \sum_{x_2,y_2\in S} \zeta(x_2)( 1-\zeta(y_2)) G_{\xi,\zeta}(x_1,y_1;x_2,y_2)\right)\nonumber\\
&&\qquad\times \bigl(f(\xi^{x_1,y_1},\zeta) - f(\xi,\zeta)\bigr)\nonumber\\
&+&\sum_{x_2,y_2 \in S} \zeta(x_2)( 1-\zeta(y_2))
\left(\Ga_\zeta(x_2,y_2)
- \sum_{x_1,y_1\in S}  \xi(x_1)( 1-\xi(y_1)) G_{\xi,\zeta}(x_1,y_1;x_2,y_2) \right) \nonumber\\
&&\qquad\times\left(f(\xi,\zeta^{x_2,y_2}) -f(\xi,\zeta)\right)\nonumber\\
&+&\sum_{x_1,y_1\in S}  \sum_{x_2,y_2\in S}  \xi(x_1)( 1-\xi(y_1))\zeta(x_2)( 1-\zeta(y_2))
G_{\xi,\zeta}(x_1,y_1;x_2,y_2)\nonumber\\
&&\qquad\times\bigl(f(\xi^{x_1,y_1},\zeta^{x_2,y_2}) - f(\xi,\zeta)\bigr)\label{6.0}
\eeq
Taking into account notations \eqref{coudis2}--\eqref{coudis3}, we get
\beq
\overline{\calL}f(\xi,\zeta)&=&
\sum_{x_1,y_1 \in S} \xi(x_1)( 1-\xi(y_1))
\bigl(\Ga_\xi(x_1,y_1)
-  \varphi_{\xi,\zeta} (x_1,y_1)\bigr)\nonumber\\
&&\qquad\times \bigl(f(\xi^{x_1,y_1},\zeta) - f(\xi,\zeta)\bigr)\nonumber\\
&+&\sum_{x_2,y_2 \in S} \zeta(x_2)( 1-\zeta(y_2))
\bigl(\Ga_\zeta(x_2,y_2)
- \overline\varphi_{\xi,\zeta} (x_2,y_2) \bigr) \nonumber\\
&&\qquad\times\bigl(f(\xi,\zeta^{x_2,y_2}) -f(\xi,\zeta)\bigr)\nonumber\\
&+&\sum_{x_1,y_1\in S}  \sum_{x_2,y_2\in S}  \xi(x_1)( 1-\xi(y_1))\zeta(x_2)( 1-\zeta(y_2))
G_{\xi,\zeta}(x_1,y_1;x_2,y_2)\nonumber\\
&&\qquad\times\bigl(f(\xi^{x_1,y_1},\zeta^{x_2,y_2}) - f(\xi,\zeta)\bigr)\label{6.1}
\eeq 
 In the above expression, the first two terms on the r.h.s.
  refer to uncoupled transitions, respectively 
$(\xi,\zeta)\rightarrow (\xi^{x_1,y_1},\zeta) $ and
$(\xi,\zeta)\rightarrow (\xi,\zeta^{x_2,y_2}) $,
while the third line refers to coupled transitions  
$(\xi,\zeta)\rightarrow (\xi^{x_1,y_1},\zeta^{x_2,y_2}) $.

Inequalities  \eqref{ineqcoupling2}--\eqref{ineqcoupling1} and non-negativity of 
$G_{\xi,\zeta}$ insure that the rates of all uncoupled and coupled transitions above are non-negative.\par
Moreover, if $f(\xi,\zeta)=g(\xi)$ depends only on $\xi$ (resp. $f(\xi,\zeta)=h(\zeta)$ 
depends only on $\zeta$),
we have $\overline{\calL}g(\xi)={\calL}g(\xi)$ (resp. $\overline{\calL}h(\zeta)={\calL}h(\zeta)$). 
Therefore $\overline{\calL}$ defines indeed a coupling of two copies of a generalized exclusion process. 
\epr

\bpr[Proposition \ref{propincreasing}]

\noindent
$\bullet$ We first consider the limits of the series defined in \eqref{SYx}--\eqref{bSXy}. 
By construction, these series are nonnegative, nondecresing and by \eqref{existence} 
there are also bounded from above. They are thus (absolutely) convergent and  
we denote their limits as
\beq\label{SX*}
S_{\xi,\zeta}^{x,*}=\lim_{n\to\infty} S_{\xi,\zeta}^{x,n} \\
\label{bTX*}
\overline T_{\xi,\zeta}^{x,*}=\lim_{n\to\infty} \overline T_{\xi,\zeta}^{x,n}\\
T_{\xi,\zeta}^{y,*} =\lim_{n\to\infty} T_{\xi,\zeta}^{y,n} 
\label{TY*}\\
\overline S_{\xi,\zeta}^{y,*} =\lim_{n\to\infty} \overline S_{\xi,\zeta}^{y,n} 
\label{bSY*}
\eeq
In these notations, equations \eqref{eq:1}--\eqref{eq:2} read:

For any configurations $\xi$,$\zeta$ in $\Om$ such that $\xi\le\zeta$,
\beq\label{eq:1ST}
\text{For all $x\in S$  such that } \xi(x)=1 , \qquad \overline T_{\xi,\zeta}^{x,*}\le S_{\xi,\zeta}^{x,*}\\
\text{For all $y\in S$  such that } \zeta(y)=0 , \qquad T_{\xi,\zeta}^{y,*} \le \overline S_{\xi,\zeta}^{y,*}
\label{eq:2ST}
\eeq
\noindent
$\bullet$ We now prove that  for any two nondecreasing, convergent
series $(S_n)_{n\ge 0}$ and   
$(T_n)_{n\ge 0}$, the quantity defined in \eqref{Hmn} $H_{m,n}(S_.,T_.)$ is nonnegative
 for all $m,n>0$. We have
\beq H_{m,n}(S_.,T_.) &=& S_m \wedge T_n - S_{m-1} \wedge T_n - S_m \wedge T_{n-1}
 + S_{m-1} \wedge T_{n-1}\nonumber\\
&=& \left(S_m \wedge T_n - S_{m-1} \wedge\left( S_m\wedge T_n\right)\right)
 - \left(S_m \wedge T_{n-1} - S_{m-1} \wedge\left( S_m\wedge \T_{n-1}\right)\right)\nonumber\\
&=& \left[S_m \wedge T_n - S_{m-1}\right]^+ - \left[S_m \wedge T_{n-1}
- S_{m-1}\right]^+\nonumber\\
&\ge& 0
\eeq
In the equations above, we used  $S_m\ge S_{m-1}$ to get the second line,
 the third line is an identity and, finally,  positivity comes from the fact that 
 $T_n\ge T_{m-1}$ and the function $t\to  \left[S_m \wedge t - S_{m-1}\right]^+$ is not decreasing.

In addition, we get that the sums $\displaystyle \sum_{m>0} H_{m,n}(S_.,T_.) $ 
and $\displaystyle \sum_{m>0} H_{n,m}(S_.,T_.) $ are absolutely convergent 
for all $n >0$ whenever the two series converge. In particular, 
setting $S_* ={\displaystyle \lim_{m\to\infty} S_m}$ and  
$T_* ={\displaystyle \lim_{n\to\infty} T_n}$, one gets
\beq\label{Hmn-s1}
\sum_{m>0}  H_{m,n}(S_.,T_.) = S_*\wedge T_n -S_*\wedge T_{n-1} \qquad\text{ for all } n>0\\
\sum_{n>0}  H_{m,n}(S_.,T_.) = S_m\wedge T_* -S_{m-1}\wedge T_* \qquad\text{ for all } m>0\label{Hmn-s2}
\eeq
We are now ready to turn to the proof of Proposition \ref{propincreasing}.  

\noindent
$\bullet$ We first prove that the coupling rates \eqref{G-rates} satisfy conditions
 \eqref{ineqcoupling1}--\eqref{ineqcoupling2} of Proposition \ref{propcoupling0}.

First, for all non ordered pairs of configurations $(\xi,\zeta)\in \Om\times\Om$, 
all coupling rates $G_{\xi,\zeta}$
defined by \eqref{G-rates} are zero, so that the left hand sides of equations 
\eqref{ineqcoupling1}--\eqref{ineqcoupling2} are identically zero and both equations
 \eqref{ineqcoupling1}--\eqref{ineqcoupling2} trivially hold.

We now consider the case $(\xi,\zeta)\in \Om\times\Om$ with $\xi\le\zeta$. 
For all $(x,y)\in S^2$, the left hand side 
of equation \eqref{ineqcoupling1} reads:
\beq
\varphi_{\xi,\zeta}(x,y)
&=&
\sum_{x',y'\in S} \zeta(x') \left(1-\zeta(y')\right) G_{\xi,\zeta}(x,y;x',y')\nonumber\\
&=& 
\zeta(x)\left(1-\zeta(y)\right)\Gamma_\xi(x,y)\wedge\Gamma_\zeta(x,y)\nonumber\\
&&+ \zeta(x) {\displaystyle \sum_{m,n>0} 
\delta(y,y_{\xi,\zeta}^{x,m}) \, \sum_{y'\in S} \delta(y',\overline y_{\xi,\zeta}^{x,n})} \,
 H_{m,n} (S_{\xi,\zeta}^{x,.},\overline T_{\xi,\zeta}^{x,.})\nonumber\\
&&+ \left(1-\zeta(y)\right) {\displaystyle\sum_{m,n>0}
 \delta(x,x_{\xi,\zeta}^{y,m}) \sum_{x'\in S} \delta(x',\overline x_{\xi,\zeta}^{y,n}) }\,
H_{m,n} (T_{\xi,\zeta}^{y,.},\overline S_{\xi,\zeta}^{x,.})\nonumber\\
&=& 
\zeta(x)\left(1-\zeta(y)\right)\Gamma_\xi(x,y)\wedge\Gamma_\zeta(x,y)\nonumber\\
&&+ \zeta(x)\zeta(y) {\displaystyle \sum_{m>0} }
\delta(y,y_{\xi,\zeta}^{x,m}) \, 
\left(S_{\xi,\zeta}^{x,m}\wedge \overline T_{\xi,\zeta}^{x,*} 
- S_{\xi,\zeta}^{x,m-1}\wedge \overline T_{\xi,\zeta}^{x,*}\right)
 \nonumber\\
&&+ \zeta(x)\left(1-\zeta(y)\right) {\displaystyle\sum_{m>0}}
 \delta(x,x_{\xi,\zeta}^{y,m}) 
\left(T_{\xi,\zeta}^{y,m}\wedge\overline S_{\xi,\zeta}^{y,*}
 -T_{\xi,\zeta}^{y,m-1}\wedge \overline S_{\xi,\zeta}^{y,*}\right) \nonumber\\
%
%
&=& 
\zeta(x)\left(1-\zeta(y)\right)\Gamma_\xi(x,y)\wedge\Gamma_\zeta(x,y)\nonumber\\
&&+ \zeta(x)\zeta(y) {\displaystyle \sum_{m>0} }
\delta(y,y_{\xi,\zeta}^{x,m}) \, 
\left(S_{\xi,\zeta}^{x,m}\wedge \overline T_{\xi,\zeta}^{x,*} 
- S_{\xi,\zeta}^{x,m-1}\wedge \overline T_{\xi,\zeta}^{x,*}\right)
 \nonumber\\
&&+ \zeta(x)\left(1-\zeta(y)\right) {\displaystyle\sum_{m>0}}
 \delta(x,x_{\xi,\zeta}^{y,m}) 
\left[\Gamma_\xi(x,y)-\Gamma_\zeta(x,y)\right]^+ \nonumber\\
&=& 
\zeta(x)\left(1-\zeta(y)\right)\Gamma_\xi(x,y)\nonumber\\
&&+ \zeta(x)\zeta(y) {\displaystyle \sum_{m>0} }
\delta(y,y_{\xi,\zeta}^{x,m}) \, 
\left(S_{\xi,\zeta}^{x,m}\wedge \overline T_{\xi,\zeta}^{x,*} 
- S_{\xi,\zeta}^{x,m-1}\wedge \overline T_{\xi,\zeta}^{x,*}\right)
 \label{eq:marg1a}
\eeq
In the second expression, the summation over $y'$ in the second term 
and the summation over $x'$ in the third term just give 1 and we use 
the expression \eqref{Hmn-s2} to compute the summation over $n>0$. 
The fourth equality is a consequence of relation \eqref{eq:1ST}, which gives
\[
T_{\xi,\zeta}^{y,m}\wedge\overline S_{\xi,\zeta}^{y,*} 
-T_{\xi,\zeta}^{y,m-1}\wedge \overline S_{\xi,\zeta}^{y,*}
= T_{\xi,\zeta}^{y,m} - T_{\xi,\zeta}^{y,m-1}
=  \left[\Gamma_\xi(x_{\xi,\zeta}^{y,m}, y)- \Gamma_\zeta(x_{\xi,\zeta}^{y,m}, y)\right]^+
\]
Now using the estimate
\be
S_{\xi,\zeta}^{x,m}\wedge \overline T_{\xi,\zeta}^{x,*} - S_{\xi,\zeta}^{x,m-1}\wedge \overline T_{\xi,\zeta}^{x,*}
\le S_{\xi,\zeta}^{x,m} - S_{\xi,\zeta}^{x,m-1} = \Gamma_\xi(x, y_{\xi,\zeta}^{x,m})\nonumber
\ee
we get the inequality
\be
\varphi_{\xi,\zeta}(x,y)
\le 
\zeta(x)\left(1-\zeta(y)\right)\Gamma_\xi(x,y)
+ \zeta(x)\zeta(y) \Gamma_\xi(x,y)\le  \Gamma_\xi(x,y)\label{ineq:marg1}
\ee
%
%
%
%
Thus inequality  \eqref{ineqcoupling1} holds for $\xi\le\zeta$.

We prove \eqref{ineqcoupling2}  for $\xi\le\zeta$ in a similar way, as follows. 
 For all $(x,y)\in S^2$, the left hand side 
of equation \eqref{ineqcoupling2} reads:
\beq
\overline\varphi_{\xi,\zeta}(x,y)
&=&
\sum_{x',y'\in S} \xi(x') \left(1-\xi(y')\right) G_{\xi,\zeta}(x',y';x,y)\nonumber\\
&=& 
\xi(x)\left(1-\xi(y)\right)\Gamma_\xi(x,y)\wedge\Gamma_\zeta(x,y)\nonumber\\
&&+ \xi(x) {\displaystyle \sum_{m,n>0} \sum_{y'\in S}
\delta(y',y_{\xi,\zeta}^{x,m}) \,  \delta(y,\overline y_{\xi,\zeta}^{x,n})} \,
 H_{m,n} (S_{\xi,\zeta}^{x,.},\overline T_{\xi,\zeta}^{x,.})\nonumber\\
&&+ \left(1-\xi(y)\right) {\displaystyle\sum_{m,n>0} \sum_{x'\in S}
 \delta(x',x_{\xi,\zeta}^{y,m})  \delta(x,\overline x_{\xi,\zeta}^{y,n}) }\,
H_{m,n} (T_{\xi,\zeta}^{y,.},\overline S_{\xi,\zeta}^{x,.})\nonumber\\
&=& 
\xi(x)\left(1-\xi(y)\right)\Gamma_\xi(x,y)\wedge\Gamma_\zeta(x,y)\nonumber\\
&&+ \xi(x)\left(1-\xi(y)\right) {\displaystyle \sum_{n>0} }
 \delta(y,\overline y_{\xi,\zeta}^{x,n}) \, 
\left(S_{\xi,\zeta}^{x,*}\wedge \overline T_{\xi,\zeta}^{x,n} 
- S_{\xi,\zeta}^{x,*}\wedge \overline T_{\xi,\zeta}^{x,n-1}\right)
 \nonumber\\
&&+ \left(1-\xi(x)\right)\left(1-\xi(y)\right) {\displaystyle\sum_{n>0}}
  \delta(x,\overline x_{\xi,\zeta}^{y,n})
\left(T_{\xi,\zeta}^{y,*}\wedge\overline S_{\xi,\zeta}^{y,n} 
-T_{\xi,\zeta}^{y,*}\wedge \overline S_{\xi,\zeta}^{y,n-1}\right) \nonumber\\
&=& 
\xi(x)\left(1-\xi(y)\right)\Gamma_\xi(x,y)\wedge\Gamma_\zeta(x,y)\nonumber\\
&&+ \xi(x)\left(1-\xi(y)\right) {\displaystyle \sum_{n>0} }
 \delta(y,\overline y_{\xi,\zeta}^{x,n}) \, 
\left[\Gamma_\zeta(x,y)\wedge\Gamma_\xi(x,y) \right]^+
 \nonumber\\
&&+ \left(1-\xi(x)\right)\left(1-\xi(y)\right) {\displaystyle\sum_{n>0}}
  \delta(x,\overline x_{\xi,\zeta}^{y,n})
\left(T_{\xi,\zeta}^{y,*}\wedge\overline S_{\xi,\zeta}^{y,n} 
-T_{\xi,\zeta}^{y,*}\wedge \overline S_{\xi,\zeta}^{y,n-1}\right) \nonumber\\
&=& 
\xi(x)\left(1-\xi(y)\right)\Gamma_\zeta(x,y)\nonumber\\
&&+ \left(1-\xi(x)\right)\left(1-\xi(y)\right) {\displaystyle\sum_{n>0}}
  \delta(x,\overline x_{\xi,\zeta}^{y,n})
\left(T_{\xi,\zeta}^{y,*}\wedge\overline S_{\xi,\zeta}^{y,n} 
-T_{\xi,\zeta}^{y,*}\wedge \overline S_{\xi,\zeta}^{y,n-1}\right)
\eeq

In the second expression, the summation over $y'$ in the second term 
and the summation over $x'$ in the third term just give 1 and we use 
the expression \eqref{Hmn-s1} to compute the summation over $n>0$. 
To get the fourth expression, we used the relation \eqref{eq:2ST} to obtain 
\be
S_{\xi,\zeta}^{x,*}\wedge \overline T_{\xi,\zeta}^{x,n} 
- S_{\xi,\zeta}^{x,*}\wedge \overline T_{\xi,\zeta}^{x,n-1}
=\overline T_{\xi,\zeta}^{x,n} - \overline T_{\xi,\zeta}^{x,n-1}= 
\left[\Gamma_\zeta(x, \overline y_{\xi,\zeta}^{x,n})- \Gamma_\xi(x, \overline y_{\xi,\zeta}^{x,n})\right]^+
\nonumber
\ee
Now we have the estimate
\be
T_{\xi,\zeta}^{y,*}\wedge\overline S_{\xi,\zeta}^{y,n} 
-T_{\xi,\zeta}^{y,*}\wedge \overline S_{\xi,\zeta}^{y,n-1}
\le \overline S_{\xi,\zeta}^{y,n} - \overline S_{\xi,\zeta}^{y,n-1} 
= \Gamma_\zeta(\overline x_{\xi,\zeta}^{y,n}, y)
\nonumber
\ee
%
which gives 
\beq
\overline\varphi_{\xi,\zeta}(x,y)
&\le& 
\xi(x)\left(1-\xi(y)\right)\Gamma_\zeta(x,y)
+ \left(1-\xi(x)\right)\left(1-\xi(y)\right) \Gamma_\zeta(x,y)
\le \Gamma_\zeta(x,y)\label{ineq:marg2}
\eeq

Equation \eqref{ineqcoupling2} is proven for $\xi\le\zeta$.

A similar derivation holds in the case $\zeta<\xi$. Thus the coupling rates defined 
in Proposition \ref{propincreasing} satisfy the conditions  
\eqref{ineqcoupling1}--\eqref{ineqcoupling2} of Proposition \ref{propcoupling0}. \\ \\

\noindent
$\bullet$ We now prove that this coupling is increasing.

We suppose that $\xi\le\zeta$. We first consider  coupled  transitions. 
From equation \eqref{G-rates}, we find that
a coupled transition $(\xi,\zeta) \to (\xi^{x,y}, \zeta^{x',y'})$ 
has possibly a non zero coupling rate $G_{\xi,\zeta}(x,y;x',y')$ in three possible cases:
\begin{itemize}
\item $x=x'$ and $y=y'$:

thus
\beq
\xi^{x,y}(x) &=& \zeta^{x',y'}(x)=0\nonumber\\
\xi^{x,y}(y) &=& \zeta^{x',y'}(y)=1\nonumber\\
\xi^{x,y}(z) &=& \xi(z) \le \zeta(z) = \zeta^{x',y'}(z) \qquad\text{ for all } z\not= x,y\nonumber
\eeq 
\item $x=x'$, $y\in Y_{\xi,\zeta}^x$ and $y'\in \overline Y_{\xi,\zeta}^x$

thus  $y\not=y',\,\zeta(y)=1$ and
\beq
\xi^{x,y}(x) &=& \zeta^{x',y'}(x)=0\nonumber\\
\xi^{x,y}(y) &\le& 1 = \zeta (y) = \zeta^{x',y'}(y)\nonumber\\
\xi^{x,y}(y') &\le& 1 =  \zeta^{x',y'}(y')\nonumber\\
\xi^{x,y}(z) &=& \xi(z) \le \zeta(z) = \zeta^{x',y'}(z) \qquad\text{ for all } z\not= x,y, y'\nonumber
\eeq 
\item $y=y'$, $x\in X_{\xi,\zeta}^y$ and $x'\in \overline X_{\xi,\zeta}^y$

thus  $x\not=x',\,\xi(x') =0$ and
\beq
\xi^{x,y}(x) &=& 0 \le  \zeta^{x',y'}(x)\nonumber\\
\xi^{x,y}(x') &=& \xi(x') = 0 \le \zeta^{x',y'}(x')\nonumber\\
\xi^{x,y}(y) &=&\zeta^{x',y'}(y') =1\nonumber\\
\xi^{x,y}(z) &=& \xi(z) \le \zeta(z) = \zeta^{x',y'}(z) \qquad\text{ for all } z\not= x,x',y\nonumber
\eeq  
\end{itemize}
In all three cases, we find that $\xi^{x,y}\le \zeta^{x',y'}$. 
Hence partial order is preserved in coupled transitions 
for all $\xi\le\zeta$.

We now turn to uncoupled transitions, $(\xi,\zeta) \to (\xi^{x,y}, \zeta)$ 
and $(\xi,\zeta) \to (\xi, \zeta^{x',y'})$, with rates 
$\left(\Gamma_\xi (x,y) -\varphi_{\xi,\zeta} (x,y)\right)$ and 
$\left(\Gamma_\zeta (x,y) -\overline \varphi_{\xi,\zeta} (x,y)\right)$
respectively. In both cases,  partial order could be broken if and only if 
$\xi(x) =\zeta(x) =1$, $\xi(y)=\zeta(y) =0$ and the associated transition rate 
is nonzero. In the first case, $\zeta(x)=1$ implies that  $y\not\in Y_{\xi,\zeta}^x$, 
which allows us to  precise the estimate  \eqref{ineq:marg1} and get the value of 
$\varphi_{\xi,\zeta} (x,y)$, as  follows. Note that 
in the expression \eqref{G-rates} for $G _{\xi,\zeta}(x,y;x',y')$ when $\xi\leq\zeta$,
since  $y\not\in Y_{\xi,\zeta}^x$,
we are in the case $y=y'$ so that $G _{\xi,\zeta}(x,y;x',y')$ is given by the third line
in \eqref{G-rates}: 
\beq
\varphi_{\xi,\zeta}(x,y)
&=&
\Gamma_\xi(x,y)\wedge\Gamma_\zeta(x,y)
+  {\displaystyle\sum_{m,n>0}
 \delta(x,x_{\xi,\zeta}^{y,m}) \sum_{x'\in S} \delta(x',\overline x_{\xi,\zeta}^{y,n}) }\,
H_{m,n} (T_{\xi,\zeta}^{y,.},\overline S_{\xi,\zeta}^{x,.})\nonumber\\
&=& 
\Gamma_\xi(x,y)\wedge\Gamma_\zeta(x,y)+  {\displaystyle\sum_{m>0}}
 \delta(x,x_{\xi,\zeta}^{y,m}) 
\left(T_{\xi,\zeta}^{y,m}\wedge\overline S_{\xi,\zeta}^{y,*} 
-T_{\xi,\zeta}^{y,m-1}\wedge \overline S_{\xi,\zeta}^{y,*}\right) \nonumber\\
&=& 
\Gamma_\xi(x,y)\wedge\Gamma_\zeta(x,y)+ {\displaystyle\sum_{m>0}}
 \delta(x,x_{\xi,\zeta}^{y,m}) 
\left[\Gamma_\xi(x,y) -\Gamma_\zeta(x,y)\right]^+ \nonumber\\
\nonumber
&=& 
\Gamma_\xi(x,y)\wedge\Gamma_\zeta(x,y)
+  {\bf 1}_{x\in X_{\xi,\zeta}^y}
\left[\Gamma_\xi(x,y) -\Gamma_\zeta(x,y)\right]^+ \nonumber\\
&=& 
\Gamma_\xi(x,y)\wedge\Gamma_\zeta(x,y)+ 
\left[\Gamma_\xi(x,y) -\Gamma_\zeta(x,y)\right]^+ \nonumber\\
&=& \Gamma_\xi(x,y)\label{eq:marg1}
\eeq 
Thus uncoupled transitions in the first marginal that do not preserve partial order 
in the case $\xi\le\zeta$ have zero transition rates.

In uncoupled transitions for the second marginal, $(\xi,\zeta) \to (\xi, \zeta^{x,y})$ 
in which partial order could be broken, $\xi(y)=0$ implies $x\not\in X_{\xi,\zeta}^y$ and, 
following the same line as in \eqref{ineq:marg2}, one gets now the value of 
$\overline \varphi_{\xi,\zeta} (x,y)$.  Again,  in the expression \eqref{G-rates} 
for $G _{\xi,\zeta}(x,y;x',y')$ when $\xi\leq\zeta$,
since  $x\not\in X_{\xi,\zeta}^y$,
we are in the case $x=x'$ so that $G _{\xi,\zeta}(x,y;x',y')$ is given by the second line
in \eqref{G-rates}: 
\beq
\overline\varphi_{\xi,\zeta}(x,y)
&=&
\sum_{x',y'\in S} \xi(x') \left(1-\xi(y')\right) G_{\xi,\zeta}(x',y';x,y)\nonumber\\
&=& 
\Gamma_\xi(x,y)\wedge\Gamma_\zeta(x,y
+  {\displaystyle \sum_{m,n>0} \sum_{y'\in S}
\delta(y',y_{\xi,\zeta}^{x,m}) \,  \delta(y,\overline y_{\xi,\zeta}^{x,n})} \,
 H_{m,n} (S_{\xi,\zeta}^{x,.},\overline T_{\xi,\zeta}^{x,.})\nonumber\\
&=& 
\Gamma_\xi(x,y)\wedge\Gamma_\zeta(x,y)+ {\displaystyle \sum_{n>0} }
 \delta(y,\overline y_{\xi,\zeta}^{x,n}) \, 
\left(S_{\xi,\zeta}^{x,*}\wedge \overline T_{\xi,\zeta}^{x,n} 
- S_{\xi,\zeta}^{x,*}\wedge \overline T_{\xi,\zeta}^{x,n-1}\right)
 \nonumber\\
&=& 
\Gamma_\xi(x,y)\wedge\Gamma_\zeta(x,y)+{\displaystyle \sum_{n>0} }
\delta(y,\overline y_{\xi,\zeta}^{x,n}) \, 
\left[\Gamma_\zeta(x,y) -\Gamma_\xi(x,y)\right]^+
 \nonumber\\
&=& 
\Gamma_\xi(x,y)\wedge\Gamma_\zeta(x,y)+ {\bf 1}_{y\in \overline Y_{\xi,\zeta}^x} 
\left[\Gamma_\zeta(x,y) -\Gamma_\xi(x,y)\right]^+ \nonumber\\
&=& 
\Gamma_\xi(x,y)\wedge\Gamma_\zeta(x,y)+\,\left[\Gamma_\zeta(x,y) -\Gamma_\xi(x,y)\right]^+ 
 \nonumber\\
&=& \Gamma_\zeta(x,y)\label{eq:marg2}
\eeq 
Uncoupled transitions in the second marginal that do not preserve partial order 
in the case $\xi\le\zeta$ have thus also zero transition rates.

In conclusion, in the generator of the coupling process \eqref{coupling}
 with rates \eqref{G-rates},  for all pairs of configurations $(\xi,\zeta)\in \Om\times \Om$ 
 such that $\xi\le\zeta$,  all possible transitions, coupled or uncoupled which have a non zero 
 transition rate do preserve the partial order. In the case $\xi> \zeta$, the same result
can be obtained along similar lines, and we thus omit its proof. The coupling defined 
in Proposition  \ref{propincreasing} is thus increasing.
\epr

\bpr[Proposition \ref{propdecdisc}]

\noindent$\bullet$ We first prove that the operator defined by \eqref{gene2} is a valid coupling, 
that is the coefficient associated to each transition is nonnegative.
We rewrite the generator  $\overline{\calL}^D$  as,
\beq\label{pr:gene2}
\overline{\calL}^Df(\xi,\zeta)&=&
\sum_{x_1,y_1 \in S} \xi(x_1)( 1-\xi(y_1))
\bigl(\Ga_\xi(x_1,y_1)
-  \sum_{x_2,y_2\in S} \zeta(x_2)( 1-\zeta(y_2)) G^D_{\xi,\zeta}(x_1,y_1;x_2,y_2)\bigr)\nonumber\\
&&\qquad\times \bigl(f(\xi^{x_1,y_1},\zeta) - f(\xi,\zeta)\bigr)\nonumber\\
&+&\sum_{x_2,y_2 \in S} \zeta(x_2)( 1-\zeta(y_2))
\bigl(\Ga_\zeta(x_2,y_2)
- \sum_{x_1,y_1\in S}  \xi(x_1)( 1-\xi(y_1)) G^D_{\xi,\zeta}(x_1,y_1;x_2,y_2) \bigr) \nonumber\\
&&\qquad\times\bigl(f(\xi,\zeta^{x_2,y_2}) -f(\xi,\zeta)\bigr)\\
&+&\sum_{x_1,y_1\in S}  \sum_{x_2,y_2\in S}  \xi(x_1)( 1-\xi(y_1))\zeta(x_2)( 1-\zeta(y_2))
G^D_{\xi,\zeta}(x_1,y_1;x_2,y_2)\nonumber\\
&&\qquad\times\bigl(f(\xi^{x_1,y_1},\zeta^{x_2,y_2}) - f(\xi,\zeta)\bigr)\nonumber
\eeq
In the above expression, the first (respectively second) line refers to uncoupled transitions
$(\xi,\zeta)\rightarrow (\xi^{x_1,y_1},\zeta) $ (respectively 
$(\xi,\zeta)\rightarrow (\xi,\zeta^{x_2,y_2}) $),
while the third line refers to coupled transitions  
$(\xi,\zeta)\rightarrow (\xi^{x_1,y_1},\zeta^{x_2,y_2}) $.

We first prove that the coefficient associated to an uncoupled transition 
$(\xi,\zeta)\rightarrow (\xi^{x_1,y_1},\zeta) $ is non-negative. It reads
\beq
&&\xi(x_1)( 1-\xi(y_1))
\bigl(\Ga_\xi(x_1,y_1)
-  \sum_{x_2,y_2\in S} \zeta(x_2)( 1-\zeta(y_2)) G^D_{\xi,\zeta}(x_1,y_1;x_2,y_2)\bigr)\nonumber\\
&&= \xi(x_1)( 1-\xi(y_1))
\bigl(\Ga_\xi(x_1,y_1)
-  \sum_{x_2,y_2\in S} \zeta(x_2)( 1-\zeta(y_2))
\nonumber\\
&&\qquad\times  \sum_{x,y\in S}(\xi\vee\zeta)(x) (1-(\xi\vee\zeta)(y))
 {\displaystyle \frac{1}{N_{\xi,\zeta}(x,y)}} G_{\xi,\xi\vee\zeta}(x_1,y_1;x,y) G_{\xi\vee\zeta,\zeta}(x,y;x_2,y_2)\bigr)\\
&&=  \xi(x_1)( 1-\xi(y_1)) \bigl(\Ga_\xi(x_1,y_1)\nonumber\\
&&\qquad -  \sum_{x,y\in S}(\xi\vee\zeta)(x) (1- (\xi\vee\zeta)(y))
  {\displaystyle \frac{1}{N_{\xi,\zeta}(x,y)}} \varphi_{\xi\vee\zeta,\zeta}(x,y)\;G_{\xi,\xi\vee\zeta}(x_1,y_1;x,y) \bigr)\nonumber\\
&&\geq  \xi(x_1)( 1-\xi(y_1))
\bigl(\Ga_\xi(x_1,y_1)
-  \sum_{x,y\in S}(\xi\vee\zeta)(x) (1-(\xi\vee\zeta)(y))  G_{\xi,\xi\vee\zeta}(x_1,y_1;x,y)\bigr)\nonumber\\
&&\geq 0\nonumber
 \eeq
 In this derivation,  we used \eqref{GD} to get the first equality, 
then exchanged the summations and used \eqref{coudis3}
  to get the second one; first inequality comes from 
  $ {\displaystyle \frac{1}{N_{\xi,\zeta}(x,y)}}\, \varphi_{\xi\vee\zeta,\zeta}(x,y) \le 1$ 
  and nonnegativity of the coupling rates $G_{\xi,\xi\vee\zeta}$; 
  the last one follows from 
inequality \eqref{ineqcoupling1}. Non negativity of the coefficients 
associated to uncoupled transitions
$(\xi,\zeta)\rightarrow (\xi,\zeta^{x_2,y_2}) $ follows along similar 
lines and inequality \eqref{ineqcoupling2}. Non-negativity of $G_{\xi,\zeta}$ 
insures that the rates $G^D_{\xi,\zeta}$ of  coupled transitions 
$(\xi,\zeta)\rightarrow (\xi^{x_1,y_1},\zeta^{x_2,y_2}) $  are also non negative.

\noindent$\bullet$ We now prove that the new coupling is increasing. 

Suppose that $\xi \le\zeta$. We have $\xi\vee\zeta=\zeta$; equations \eqref{coudis2}--\eqref{coudis3} 
and Remark \ref{xi=zeta} give
\beq
\overline \varphi_{\xi,\xi\vee\zeta}(x,y)&=&\overline\varphi_{\xi,\zeta}(x,y)=
\sum_{x',y'}\xi(x')(1-\xi(y')) G_{\xi,\zeta}(x',y';x,y)\nonumber\\
\varphi_{\xi\vee\zeta,\zeta}(x,y)&=&\varphi_{\zeta,\zeta}(x,y)=
\sum_{x',y'}\zeta(x')(1-\zeta(y')) G_{\zeta,\zeta}(x,y;x',y')= \zeta(x)(1-\zeta(y))  \Ga_\zeta(x,y)\nonumber
\eeq
Inequality  \eqref{ineqcoupling2} implies here  that $\overline \varphi_{\xi,\zeta}(x,y) \le\Ga_\zeta(x,y)$,
and we get from equation \eqref{coudis1} 
\beq\label{6.27}
\zeta(x) (1-\zeta(y))\,  {\displaystyle \frac{1}{N_{\xi,\zeta}(x,y)}}\,  \Ga_{\zeta}(x,y)=
\begin{cases}
1 
& \hbox{ if } \zeta(x)(1-\zeta(y)) \Ga_\zeta(x,y) \not= 0\cr
0 & \hbox { otherwise }
\end{cases}
\eeq
This enables us to prove Remark  \ref{rk:GD=G}. For all $(x_1,y_1,x_2,y_2)\in S^4$
such that $\xi(x_1)(1-\xi(y_1))\zeta(x_2)(1-\zeta(y_2))\not=0$, the coupling rates 
$G^D_{\xi,\zeta}(x_1,y_1;x_2,y_2) $ thus read
\beq\label{6.28}
G^D_{\xi,\zeta} (x_1,y_1;x_2,y_2) &= &
\sum_{x,y\in S} \zeta(x)(1-\zeta(y)) \, {\displaystyle \frac{1}{N_{\xi,\zeta}(x,y)}}\, G_{\xi,\zeta}(x_1,y_1;x,y) 
G_{\zeta,\zeta}(x,y;x_2,y_2)
\nonumber\\
&=& \zeta(x_2)(1-\zeta(y_2))  {\displaystyle \frac{1}{N_{\xi,\zeta}(x_2,y_2)}} \Ga_{\zeta}(x_2,y_2) G_{\xi,\zeta}(x_1,y_1;x_2,y_2)
\nonumber\\
&=& 
 G_{\xi,\zeta}(x_1,y_1;x_2,y_2) 
\eeq
First equality is Equation \eqref{GD} in the present case; second equality follows from Remark \ref{xi=zeta};  
the last one follows from \eqref{6.27} and $0\le G_{\xi,\zeta}(x_1,y_1;x_2,y_2)\le  \Ga_{\zeta}(x_2,y_2) $
 (this last inequality comes from the fact that in $\overline\calL$,  the increasing coupling generator
defined in Proposition  \ref{propincreasing}, the rates of uncoupled transitions are non-negative, 
cf. \eqref{6.0}--\eqref{6.1}). 

 Inserting Equation \eqref{6.28} in \eqref{gene2}, we get
\beq
{\overline \calL}^D f(\xi,\zeta)&=& \sum_{x_1,y_1 \in S}\xi(x_1) (1-\xi(y_1))
\Ga_\xi(x_1,y_1)\bigl(f(\xi^{x_1,y_1},\zeta) - f(\xi,\zeta)\bigr)\nonumber\\
&&+\sum_{x_2,y_2 \in S}\zeta(x_2) (1-\zeta(y_2)) 
\Ga_\zeta(x_2,y_2)\bigl(f(\xi,\zeta^{x_2,y_2}) - f(\xi,\zeta)\bigr)\nonumber\\
&&+\sum_{x_1,y_1\in S}  \sum_{x_2,y_2\in S} 
\xi(x_1) (1-\xi(y_1)) \zeta(x_2) (1-\zeta(y_2)) G_{\xi,\zeta}(x_1,y_1;x_2,y_2)\nonumber\\
&&\qquad\times \bigl( f(\xi^{x_1,y_1},\zeta^{x_2,y_2})- f(\xi^{x_1,y_1},\zeta) 
- f(\xi,\zeta^{x_2,y_2})+ f(\xi,\zeta)\bigr)\nonumber\\
&=& {\overline \calL} f(\xi,\zeta)\nonumber
\eeq
%
A similar identity holds for $\xi>\zeta$. Since both generators identify on 
$\{\xi\le\zeta\}\cup \{\xi >\zeta\}$, the coupling with generator 
${\overline \calL}^D$ is also increasing.

\noindent$\bullet$ We now prove that discrepancies cannot increase 
under ${\overline \calL}^D$.

For any finite domain $D\subset S$, the number of discrepancies in $D$ 
between two configurations $\xi$, $\zeta$ in $\Om$ is defined
as
\beq
\sum_{x\in D}|\xi(x)-\zeta(x)|\nonumber
\eeq
Each transition in  \eqref{gene2} with positive transition rate involves a change 
on a finite number of sites. For any such transition, say
$(\xi,\zeta)\longrightarrow (\xi',\zeta')$, and for any finite domain $D$ 
which contains 
all sites involved in the transition
\beq\label{domdisc}
D\supset\bigl\{x\in S,\xi'(x)\not=\xi(x) \hbox{ or } \zeta'(x)\not=\zeta(x)\bigr\}
\eeq
the variation of discrepancies is
 \beq\label{vardisc2}
 \Delta_D(\xi,\zeta;\xi',\zeta')
 &=&\sum_{x\in D}|\xi'(x)-\zeta'(x)|-\sum_{x\in D}|\xi(x)
 -\zeta(x)|\nonumber\\
 &=&\sum_{x\in D}\bigl(2\xi'(x)\vee\zeta'(x)-\xi'(x)-\zeta'(x)\big)
 -\sum_{x\in D}\bigl(2\xi(x)\vee\zeta(x)-\xi(x)-\zeta(x)\big)\nonumber\\
 &=& 2\sum_{x\in D} \bigl(\xi'(x)\vee\zeta'(x)-\xi(x)\vee\zeta(x)\big)
 \eeq
 The last equality holds since the process is conservative.
 
\noindent  
  $\bullet$ We consider first a coupled transition
   $(\xi,\zeta)\longrightarrow (\xi^{x_1,y_1},\zeta^{x_2,y_2})$ 
   for some $(x_1,y_1)$, $(x_2,y_2)$ in $S^2$  with positive transition rate in $\overline\calL^D$, 
 \beq
 \xi(x_1)  (1-\xi(y_1))  \zeta(x_2)  (1-\zeta(y_2)) G^D_{\xi,\zeta}(x_1,y_1;x_2,y_2) >0
 \eeq
 Turning to the definition  \eqref{GD}, $G^D_{\xi,\zeta}(x_1,y_1;x_2,y_2) > 0$ implies  that there exists
$(x_0,y_0)\in S^2$ such that both  $G_{\xi,\xi\vee\zeta}(x_1,y_1;x_0,y_0)>0$ and 
$ G_{\xi\vee\zeta,\zeta}(x_0,y_0;x_2,y_2)>0$. Thus the transitions 
$(\xi,\xi\vee\zeta)\longrightarrow (\xi^{x_1,y_1},(\xi\vee\zeta)^{x_0,y_0})$ and  
$(\xi\vee\zeta,\zeta)\longrightarrow ((\xi\vee\zeta)^{x_0,y_0},\zeta^{x_2,y_2})$ 
have positive transition rate in $\overline \calL$. Since it is  the generator 
of an increasing coupling, $\xi\le(\xi\vee\zeta)$ and $\zeta\le(\xi\vee\zeta)$ implies that 
$\xi^{x_1,y_1} \le (\xi\vee\zeta)^{x_0,y_0}$ and $\zeta^{x_2,y_2}\le\xi\vee\zeta)^{x_0,y_0}$ and thus
\[
\xi^{x_1,y_1}\vee\zeta^{x_2,y_2}\le(\xi\vee\zeta)^{x_0,y_0}
\]
Now, for any domain $D$ as in \eqref{domdisc}, 
 \begin{eqnarray*}\label{vardisc2-2}
 \Delta_D(\xi,\zeta;\xi^{x_1,y_1},\zeta^{x_2,y_2})
 &=& \Delta_{D\cup\{x_0,y_0\}}(\xi,\zeta;\xi^{x_1,y_1},\zeta^{x_2,y_2})
\\
 &=& 2\sum_{x\in D\cup\{x_0,y_0\}} 
 \bigl(\xi^{x_1,y_1}(x)\vee\zeta^{x_2,y_2}(x)-\xi(x)\vee\zeta(x)\big)\nonumber\\
&\le& 2\sum_{x\in D\cup\{x_0,y_0\}} 
\bigl((\xi\vee\zeta)^{x_0,y_0}(x)-(\xi\vee\zeta)(x)\big)\nonumber\\
&=&0\nonumber
 \end{eqnarray*}
where the last equality follows from particle conservation.
Thus the number of discrepancies does not increase in any coupled transition 
in $\overline\calL^D$.

\noindent
 $\bullet$ We now turn to uncoupled transitions in $\overline\calL^D$. Let us consider 
 a transition in the first marginal, say $(\xi,\zeta)\longrightarrow (\xi^{x_1,y_1},\zeta)$ 
 for some $(x_1,y_1)$ in $S^2$.
For any finite domain $D$ such that $\{x_1,y_1\}\subset D$, the variation in the number
 of discrepancies reads
 \beq
 \Delta_D(\xi,\zeta;\xi^{x_1,y_1},\zeta)
 &=& 2\bigl(\xi^{x_1,y_1}(x_1)\vee\zeta(x_1)-\xi(x_1)\vee\zeta(x_1)\bigr)
 \nonumber\\&&\quad
+2\bigl(\xi^{x_1,y_1}(y_1)\vee\zeta(y_1)-\xi(y_1)\vee\zeta(y_1)\bigr)\nonumber\\
 &=& 2 \bigl(\zeta(x_1)- \zeta(y_1)\bigr)
\eeq
Thus the variation of discrepancies is non positive except in the case where both
 $\zeta(x_1)=1$ and $\zeta(y_1)=0$. We now prove that such a transition has rate $0$ 
 in  $\overline\calL^D$:

First, since $(\xi\vee\zeta)(y_1)=0$, $y_1\notin Y^{x_1}_{\xi,\xi\vee\zeta}$ and by \eqref{G-rates}, 
for any $(x,y)\in S^2$ such that $y\not= y_1$,  $G_{\xi,\xi\vee\zeta}(x_1,y_1;x,y) = 0$.  
Furthermore, $\varphi_{\xi\vee\zeta,\zeta}(x,y_1)=
\overline \varphi_{\zeta,\xi\vee\zeta,}(x,y_1)$ and since $\zeta(y_1)=0$, equation \eqref{eq:marg2} holds 
and one has
\beq
\varphi_{\xi\vee\zeta,\zeta}(x,y_1)= \Gamma_{\xi\vee\zeta} (x,y_1)
\eeq
Now the rate for the transition  $(\xi,\zeta)\longrightarrow (\xi^{x_1,y_1},\zeta)$ in $\overline\calL^D$ reads
\beq
&&\-\-\-\-\-\xi(x_1)( 1-\xi(y_1))
\bigl(\Ga_\xi(x_1,y_1)
-  \sum_{x_2,y_2\in S} \zeta(x_2)( 1-\zeta(y_2)) G^D_{\xi,\zeta}(x_1,y_1;x_2,y_2)\bigr)\nonumber\\
&=& \xi(x_1)( 1-\xi(y_1))
\bigl(\Ga_\xi(x_1,y_1)
-  \sum_{x_2,y_2\in S} \zeta(x_2)( 1-\zeta(y_2))
\nonumber\\
&&\qquad\times  \sum_{x,y\in S}(\xi\vee\zeta)(x) (1-(\xi\vee\zeta)(y))\,
  {\displaystyle \frac{1}{N_{\xi,\zeta}(x,y)}}\, G_{\xi,\xi\vee\zeta}(x_1,y_1;x,y) G_{\xi\vee\zeta,\zeta}(x,y;x_2,y_2)\bigr)\nonumber\\
&=&  \xi(x_1)( 1-\xi(y_1)) \bigl(\Ga_\xi(x_1,y_1)\nonumber\\
&&\qquad -  \sum_{x,y\in S}(\xi\vee\zeta)(x) (1- (\xi\vee\zeta)(y))\,
  {\displaystyle \frac{1}{N_{\xi,\zeta}(x,y)}}\, \varphi_{\xi\vee\zeta,\zeta}(x,y)\;G_{\xi,\xi\vee\zeta}(x_1,y_1;x,y) \bigr)\nonumber\\
&=&  \xi(x_1)( 1-\xi(y_1))
\bigl(\Ga_\xi(x_1,y_1)
-  \sum_{x,y\in S}(\xi\vee\zeta)(x) (1-(\xi\vee\zeta)(y))  G_{\xi,\xi\vee\zeta}(x_1,y_1;x,y)\bigr)\nonumber\\
&=&  \xi(x_1)( 1-\xi(y_1))
\bigl(\Ga_\xi(x_1,y_1)
-   \varphi_{\xi, \xi\vee\zeta}(x_1,y_1)  \bigr)\nonumber\\
&=& 0\nonumber
 \eeq
where the third equality comes from the fact that $ {\displaystyle \frac{1}{N_{\xi,\zeta}(x,y)}}\, \varphi_{\xi\vee\zeta,\zeta}(x,y)=1$ if $y=y_1$ and $ \Gamma_{\xi\vee\zeta} (x,y_1) > 0$, and $G_{\xi,\xi\vee\zeta}(x_1,y_1;x,y)=0$ otherwise; the last equality comes from $(\xi\vee\zeta)(x_1)=1$, $\xi\le \xi\vee\zeta$ and equation \eqref{eq:marg1}.

Thus the number of discrepancies does not increase in any uncoupled, first marginal transition in $\overline\calL^D$.

\noindent $\bullet$ Finally, we consider an uncoupled, second marginal transition 
$(\xi,\zeta)\longrightarrow (\xi,\zeta^{x_2,y_2})$ for some $(x_2,y_2)$ in $S^2$. 
Again, one proves that either the number of discrepancies does not increase, 
or has zero transition rate. The derivation is essentially identical to the previous one so we skip it.

Collecting all cases, we have shown that in any  transition in  $\overline \calL^D$ 
with nonzero transition rate, 
the number of discrepancies does not increase. The result is proven.
\epr

\bpr [Proposition \ref{decrdiscr}]
\mbox{}\\
We first construct a new increasing coupling as in Proposition 
\ref{propincreasing}, with the new rates
\beq\label{G-rates-2}
&&G _{\xi,\zeta}(x,y;x',y') =\nonumber\\
&&\qquad\begin{cases}
 \delta(x,x')\, \delta(y,y')\, \Ga_\xi(x,y)\wedge\Ga_\zeta(x,y)\\
\qquad+  \delta(x,x') {\displaystyle {\bf 1}_{y\in Y_{\xi,\zeta}^{x}}\,{\bf 1}_{y' \in \overline Y_{\xi,\zeta}^{x}}}\,
{\displaystyle \frac{1}{N_{\xi,\zeta}^{x,*}}}\,\Gamma_\xi(x,y)\,\left[  \Gamma_\zeta(x,y') -\Gamma_\xi(x,y')\right]^+
\\
\qquad+ \delta(y,y') {\displaystyle {\bf 1}_{x\in X_{\xi,\zeta}^{y}}\,{\bf 1}_{x' \in \overline X_{\xi,\zeta}^{y}}}
 { \displaystyle \frac{1}{\overline N_{\xi,\zeta}^{y,*}}}\,\left[  \Gamma_\xi(x,y) -\Gamma_\zeta(x,y)\right]^+
\,\Gamma_\zeta(x',y)
 &\text{ if } \xi\le\zeta
\cr
 \delta(x,x') \,\delta(y,y')\, \Ga_\xi(x,y)\wedge\Ga_\zeta(x,y)\\
\qquad+  \delta(x,x'){\displaystyle {\bf 1}_{y\in \overline Y_{\zeta,\xi}^{x}}\,{\bf 1}_{y' \in  Y_{\zeta,\xi}^{x}}} \, 
{\displaystyle \frac{1}{N_{\zeta,\xi}^{x,*}}}\,\left[  
\Gamma_\xi(x,y) -\Gamma_\zeta(x,y)\right]^+\,\Gamma_\zeta(x,y')
\\
\qquad+ \delta(y,y') {\displaystyle {\bf 1}_{x\in  \overline X_{\zeta,\xi}^{y}}\,{\bf 1}_{x' \in X_{\zeta,\xi}^{y}}} \,
{\displaystyle \frac{1}{\overline N_{\zeta,\xi}^{y,*}}}\,
\Gamma_\xi(x,y)\,\left[  \Gamma_\zeta(x',y) -\Gamma_\xi(x',y)\right]^+
&\text{ if } \xi>\zeta
\cr
0  &\text{ otherwise }
\end{cases}
\eeq
where 
\be
{\displaystyle {N_{\xi,\zeta}^{x,*}}}=\begin{cases}S_{\xi,\zeta}^{x,*}&\text{ if } S_{\xi,\zeta}^{x,*} >0\cr 1&\text{otherwise } \end{cases}
\ee
and similar definitions for the others normalization factors, with $S_{\xi,\zeta}^{x,*}$, $\overline S_{\xi,\zeta}^{y,*}$, $S_{\zeta,\xi}^{x,*}$ and $\overline S_{\zeta,\xi}^{y,*}$ as in Equations  \eqref{SX*}--\eqref{bSY*}.
For two configurations $\xi$ and $\zeta$ such that $\xi\le\zeta$, one can compute easily the sum of correlated jump rates associated to a jump in a given marginal. One finds, respectively
\beq
\varphi_{\xi,\zeta}(x,y)
&=&
\sum_{x',y'\in S} \zeta(x') \left(1-\zeta(y')\right) G_{\xi,\zeta}(x,y;x',y')\nonumber\\
&=& 
\zeta(x) \left(1-\zeta(y)\right)\Gamma_\xi(x,y) +\zeta(x) \,\zeta(y)  {\displaystyle \frac{\overline T_{\xi,\zeta}^{x,*}}{ S_{\xi,\zeta}^{x,*}}}\,\Gamma_\xi(x,y)
 \label{eq:marg1b}
\eeq
and 
\beq
\overline\varphi_{\xi,\zeta}(x,y)
&=&
\sum_{x',y'\in S} \xi(x') \left(1-\xi(y')\right) G_{\xi,\zeta}(x',y';x,y)\nonumber\\
&=& 
\xi(x) \left(1-\xi(y)\right)\Gamma_\zeta(x,y) +\left(1-\xi(x)\right) \,\left(1-\xi(y)\right)   {\displaystyle \frac{ T_{\xi,\zeta}^{y,*}}{ \overline S_{\xi,\zeta}^{y,*}}}\,\Gamma_\zeta(x,y)
 \label{eq:marg2b}
\eeq
Clearly coupled jump rates $G_{\xi,\zeta}(x,y;x',y')$ 
and uncoupled jump rates $\Gamma_\xi(x,y)-\varphi_{\xi,\zeta}(x,y)$,
$\Gamma_\zeta(x,y)- \overline\varphi_{\xi,\zeta}(x,y)$ are all nonnegative 
for $\xi\le\zeta$, and similarly for  $\xi>\zeta$.
Following the same lines as in the proof of Proposition \ref{propincreasing}, 
one finds that the above rates define an increasing Markovian coupling. 
Using these new rates, one can define as in Proposition \ref{propdecdisc} 
a new coupling $\overline {\mathcal L}^D$ such that the discrepancies do not increase. 
Now suppose that for a given pair of non  ordered configurations $\xi$ and $\zeta$, 
there is a discrepancy at site $x$, say $\xi(x)=1$ and $\zeta(x)=0$. 
Now the discrepancy 
can move alongside with the particle in the first marginal to any empty site $y$  
such that the edge $(x,y)$ is open 
at rate $\Gamma_\xi(x,y)>0 $, or to any fully occupied site $y$ such that the edge 
$(y,x)$ is open, alongside with a particle from the second marginal 
in the opposite direction with rate 
 ${\displaystyle\left(1-\frac{{\overline T}_{\xi,\zeta}^{x,*}}{S_{\xi,\zeta}^{x,*}}\right) \Gamma_\zeta(y,x)} >0$.
In this case,  pairs of discrepancies of opposite sign connected through an open path 
have positive probability to disappear. 
\epr

\noindent 
{\bf Acknowledgments.} 
We thank Lorenzo Bertini for useful discussions in the first stages of this work.
T. Gobron acknowledges support from the Labex CEMPI (ANR-11-LABX-0007-01).
 Part of this work has been conducted within the FP2M federation (CNRS FR 2036). 
\bibliographystyle{amsalpha}

\end{document}